\documentclass[12pt,english]{article}
\usepackage[T1]{fontenc}
\usepackage{textcomp}
\usepackage[latin9]{inputenc}
\usepackage{prettyref}
\usepackage{amstext}
\usepackage{amsthm}
\usepackage{amssymb}
\usepackage{amsmath,amssymb,amsthm,amsfonts,amsbsy,latexsym,dsfont}
\makeatletter
\numberwithin{figure}{section}
\usepackage{enumitem}
\theoremstyle{plain}
\newtheorem{thm}{\protect\theoremname}
\theoremstyle{plain}
\newtheorem{conjecture}{\protect\conjecturename}
\theoremstyle{plain}
\newtheorem{lem}[thm]{\protect\lemmaname}
\theoremstyle{plain}
\newtheorem{cor}[thm]{\protect\corollaryname}
\theoremstyle{plain}
\newtheorem{claim}[thm]{\protect\claimname}
\theoremstyle{plain}
\newtheorem{question}[conjecture]{\protect\questionname}
\theoremstyle{plain}
\newtheorem{ex}[thm]{\protect\examplename}

\usepackage{fullpage}
\usepackage{amsfonts}
\usepackage{pdfsync}
\usepackage{tikz}
\usepackage{prettyref}
\usepackage{color}

\usepackage{fullpage}

\usepackage{cite}

\DeclareMathOperator{\sh}{shift}
\DeclareMathOperator{\sw}{swap}

\usepackage[b]{esvect}

\usepackage{fancyhdr}
\usepackage{color}

\usepackage{accents,graphicx}

\makeatletter
\DeclareRobustCommand{\cev}[1]{%
  \mathpalette\do@cev{#1}%
}
\newcommand{\do@cev}[2]{%
  \fix@cev{#1}{+}%
  \reflectbox{$\m@th#1\vec{\reflectbox{$\fix@cev{#1}{-}\m@th#1#2\fix@cev{#1}{+}$}}$}%
  \fix@cev{#1}{-}%
}
\newcommand{\fix@cev}[2]{%
  \ifx#1\displaystyle
    \mkern#23mu
  \else
    \ifx#1\textstyle
      \mkern#23mu
    \else
      \ifx#1\scriptstyle
        \mkern#22mu
      \else
        \mkern#22mu
      \fi
    \fi
  \fi
}

\makeatother

\usepackage{babel}
\providecommand{\conjecturename}{Conjecture}
\providecommand{\corollaryname}{Corollary}
\providecommand{\lemmaname}{Lemma}
\providecommand{\theoremname}{Theorem}
\providecommand{\claimname}{Claim}
\providecommand{\questionname}{Question}
\providecommand{\examplename}{Example}

\def\HK{\color{blue}}

\newtheorem*{con}{Conjecture}

\def\BB{{\mathcal B}}

\def\FF{{\mathcal F}}

\def\mad{{\text{mad}}}
\def\mm{~\textrm{mod}^*~}
\def\COMMENT#1{}
\let\COMMENT=\footnote

\title{Results and Problems on Equitable  Coloring of  Graphs}

\author{
{{H. A. Kierstead}}\thanks{
\footnotesize {Arizona State University
Tempe, AZ, USA. E-mail: \texttt {kierstead@asu.edu}.
}}
\and
{{Alexandr Kostochka}}\thanks{
\footnotesize {University of Illinois at Urbana--Champaign, Urbana, IL 61801. E-mail: \texttt {kostochk@math.illinois.edu}.
 Rese arch %%% of this author
is supported in part by  NSF  Grant DMS-2153507 and by NSF RTG Grant DMS-1937241.
}}
\and
{{Zimu Xiang}}\thanks{University of Illinois at Urbana--Champaign, Urbana, IL 61801. E-mail: {\tt zimux2@illinois.edu}. }}

\date{\today}

\begin{document}

\maketitle

\begin{abstract}
 A proper coloring of vertices of a  graph is \emph{equitable} if the sizes of any two color classes differ by at most $1$. 
 Such colorings have many applications and are interesting by themselves.
 In this paper, we discuss the state of art and unsolved problems on equitable coloring and its list versions.
\end{abstract}



\section{Introduction}

Consider the following idealized problem. A service company employs a set $W$ of
$k$ equivalent workers. Each day the company must complete a set $V$ of $n$ orders.
The parameters $W$ and $V$ may vary from day-to-day. For various reasons (confidentiality,
custom, location, etc.) some pairs of orders should not or cannot be handled by the
same worker. In this case the pairs are \emph{incompatible}. Let $G=(V,E)$ be the
\emph{incompatibility} graph on $V$, where $xy\in E$ iff the orders $x$ and $y$
are incompatible. The company would like to assign orders to workers so that no two
incompatible orders are assigned to the same worker. In graph theoretic terms, this
assignment should be a \emph{proper coloring} $f:V\rightarrow W$. But the company
may have the additional requirement that the numbers of orders assigned to any two
workers should differ by at most one. In other words, the proper coloring $f$ should
be \emph{equitable}. Suppose $|V|=n$ and $|W|=k$.

Now the problem is whether $G$ has an equitable $k$-coloring, that is, a (proper)
$k$-coloring $f$ such that $\lfloor\frac{n}{k}\rfloor\leq|f^{-1}(w)|\leq\lceil\frac{n}{k}\rceil=:m$
for all colors/workers $w\in W$. Notice that this is a decision problem, not an
optimization problem; the company already employs $k$ workers. Being equitable means
that any two of the $k$ workers handle the same number of jobs $(\pm1)$. However
there are other interesting questions. For instance, suppose $m$ is a limit on the
number of orders that a worker can handle. Then on days that $n\geq km$ the company
might hire additional workers. If $G$ has an equitable $k$-coloring, does it have
an equitable $k'$-coloring for $k'>k$? Maybe the company employs a fixed set of
$k=4$ workers, but takes as long as is needed to complete all orders. When can orders
be equitably assigned to these four workers?

This example hints on possible
 applications of equitable coloring and similar notions in graph theory. Among such applications are scheduling in communication systems, construction
timetables, mutual exclusion
scheduling problem, and round-a-clock scheduling; see, e.g.~\cite{BEPSW,HHK,IL,KW,Mey,SBG,Tu}. For example, in~\cite{Tu} equitable coloring modeled scheduling weekly routes of garbage trucks.

This concept is also  useful in studying extremal combinatorial and
probabilistic problems. Alon and F\" uredi~\cite{AF} used it to study existence of some spanning subgraphs in random graphs. Alon and Yuster~\cite{AY} applied results on equitable coloring to the problem of the existence of $H$-factors in dense graphs. Janson and Ruci\' nski~\cite{JR}, Pemmaraju~\cite{Pe} and Janson~\cite{JS} used equitable colorings to
derive deviation bounds for sums of dependent random variables with limited dependency. R\" odl and Ruci\' nski~\cite{RR}  used results on equitable coloring to give a simpler proof of the Blow-Up Lemma.

\subsection{Framing the theory}

Here are two informative negative examples.

Let $S$ be  the star on $n$ vertices. It is $2$-colorable, but the maximum size of a color
class containing its root is one. Now  by the definition of an equitable coloring, every color class
has size at most two. Thus, for $S$ to have an equitable $k$-coloring, we need
$k>n/2>\Delta(S)/2$. More generally, if $v$ is a vertex of a graph $G$ on $n$
vertices such that the maximum size of an independent set containing $v$ is $m$
then every equitable coloring of $G$ uses more than $\frac{n}{m+1}$ colors. 

Let $K=K_{t,t}$. Now $K$ has an equitable $2$-coloring. If $K$ has an equitable
$t$-coloring, then every class  must have %has 
 size two and be contained in one of the two parts
of $K$. Thus $t$ is even. So when $t$ is odd, increasing the number of required
colors from $2$ to $t$ makes it impossible to equitably color $K$.

Early positive results followed from the complementary version of the problem. Let
$\overline{G}$ be the complement of $G$. Now a proper coloring of $G$ is a partition
of $V(\overline{G})$ into cliques. Suppose for simplicity that $n=km$. Then
\[
\Delta(G)<k\leftrightarrow\Delta(G)+1\leq k\leftrightarrow\delta(\overline{G})\geq n-\Delta(G)-1=(1-\frac{1}{m})n.
\]
If $m=2$ and $\Delta(G)<k$ then $\delta(\overline{G})\geq\frac{n}{2}$, so Dirac's
Theorem \cite{Dirac52} implies that $V(\overline{G})$ has a $k$-partition into
$2$-cliques. Thus $G$ has an equitable $k$-coloring. Corradi and Hajnal \cite{CH}
proved the case $m=3$: if $\delta(\overline{G})=\frac{2}{3}n$ then $V(\overline{G})$
has a $k$-partition into $3$-cliques. Then Erd\H{o}s conjectured:
\begin{con}[{Erd\H{o}s \cite{E}, 1964}]
If $H$ is a graph on $n=km$ vertices with $\delta(H)\geq n(1-\frac{1}{m})$ then
$V(H)$ has a $k$-partition into $m$-cliques.\textbf{ }
\end{con}

Hajnal offered 1000 cups of coffee for a proof of the conjecture. Gr\" unbaum \cite{Grun}
proved the case $m=4$, and reformulated the conjecture as a statement about equitable
coloring.
\begin{con}[{{Gr\" unbaum \cite{Grun}, 1968}}]
\label{conj:E}Every graph $G$ with maximum degree $\Delta(G)<k$ has a $k$-equitable
coloring.\textbf{ }
\end{con}

The Gr\" unbaum's version includes the case that $k$ does not divide $n$; this case
was handled by Larman~\cite{Larman69} (See the second paragraph of Subsection 2.1). By the
end of the sixties, the state of the art was as follows. 
\begin{thm}
\label{thm:PreHS}A graph $G$ on $n$ vertices with maximum degree $\Delta$ is
equitably k-colorable if $\Delta<k$ and either 
\begin{enumerate}[label=(\alph*),nosep]
\item \emph{(Dirac \cite{Dirac52}, 1952)}~$n/2\leq k$,
\item \emph{(Corradi and Hajnal \cite{CH}, 1963)}~$n/3\leq k$, or
\item \emph{(Gr\" unbaum \cite{Grun}, 1968)}~$n/4\leq k$.
\item \emph{(Zelinka \cite{Zelinka}, 1966)}~$k\leq3$.
\end{enumerate}
\end{thm}

Finally, in 1970 Hajnal and Szemer\' edi proved Erd\H{o}s' Conjecture (in Gr\" unbaum's
form).
\begin{thm}[Hajnal and Szemer\' edi \cite{HSz}, 1970]
\label{thm:H-Sz}Every graph with $\Delta<k$ is equitably $k$-colorable.
\end{thm}

 In fact, most of the mathematical applications of equitable coloring that we know use the  Hajnal-Szemer\' edi Theorem (mostly in the complementary form).

\medskip

Another line of research considers equitable $k$-coloring of sparse graphs $G$ with
$k<\Delta(G)$. Early results in this area include  the theorem of Meyer \cite{Mey} from 1973 on equitable 
$k$-coloring of trees $T$
for $k>\Delta(T)/2$  and the following result of Bollob\' as and Guy \cite{BG83}
on {  equitably $3$-coloring trees.}

\begin{thm}[Bollob\' as and Guy \cite{BG83}, 1983] \label{BG}
Let $T$ be a tree on $n$ vertices. Then $T$ has an equitable $k$-coloring if
 $n\geq|V(T)|\ge3\Delta(T)-8$ or
$n=3\Delta(T)-10$.
\end{thm}

Results in this spirit are discussed in Section~\ref{fewer}.

\subsection{Organization }

In  Section 2 our focus is on results for equitable $k$-coloring graphs $G$ with sufficiently
% large $k$ in terms of maximum degree. We provide a setup for proving results related
to the Hajnal-Szemer\' edi Theorem, give a complete proof of the theorem, and show that
there is a polynomial time algorithm for constructing the coloring. We generalize
Theorem~\ref{thm:H-Sz} by showing that if the Ore-degree $\Theta(G)<2k$, where
\[
\Theta(G):=\max\{d(x)+d(y):xy\in E(G)\},
\]
then $G$ has an equitable $k$-coloring, and conjecture that there is a polynomial
algorithm for finding such a coloring. Finally we discuss a strengthening of Theorem~\ref{thm:H-Sz}
to directed graphs. 

In Section 3 we discuss the Chen-Lih-Wu Conjecture and partial results 
 towards it. The
conjecture is wide open, but is proved for some interesting classes of graphs. We also discuss an equivalent conjecture and some harder conjectures, including the Ore-version of the conjecture.

%se results show that the conjecture holds for including $k\leq4$
%and $m\leq4$; further harder conjectures (hoping for counterexamples);

In Section 4 we consider problems on equitable  list coloring including
 list versions of the Hajnal-Szemer\' edi Theorem and Chen-Lih-Wu Conjecture. We also introduce  strongly equitable (SE) list coloring and compare results and problems on equitable  list coloring 
and SE list coloring. 

In Section 5  we discuss equitable coloring (and list equitable coloring) of "sparse" graphs, when we need significantly fewer colors than the maximum degree of the graph.  Finally, in Section 6 we list known results on 
some variations of equitable coloring.

\subsection{Notation}
We view
 a graph $G$ as a pair $(V,E)$ where
$V=V(G)$ denotes the set of its vertices and $E=E(G)$ denotes the set of its edges. We let $|G|=|V(G)|$ and $\|G\|=|E(G)|$. 
Further, $\Delta(G)$ is the maximum degree of $G$, and $\delta(G)$ is the minimum degree of $G$.

 For disjoint $X,Y\subseteq V(G)$,  $G[X]$ denotes the subgraph of $G$ induced by $X$,
$E(X,Y)$ denotes the set of edges in $G$ connecting $X$ with $Y$ and $\|X,Y\|=\|X,Y\|_G:=|E(X,Y)|$.
  If $X\cap Y\neq \emptyset$, then the edges of $G[X\cap Y]$ belong to $E(X,Y)$ and count twice in $\|X,Y\|.$
 If $X=\{x\}$ we may simply write $\|x,Y\|$ for $\|X,Y\|$. The degree of $x$ in $G$ is noted by $d(x)=d_{G}(x):=\|x,V(G)-x\|_{G}$.

%$G=(V,E)$. $|G|$, $\|G\|$, $G[X]$ $\widetilde{v}$, $\widetilde{\alpha}$, $E(X,Y)$,
%$E(x,Y)$ $E(X)$, $\|X,Y\|,$ including double counting when $X\cap Y\ne\emptyset$, $\|x,Y\|$, $\|X\|$, etc. 

We use Greek letters for colors and Roman letters for vertices. For a coloring $f:V\rightarrow C$
and $\gamma\in C$, set $\widetilde{\gamma}:=f^{-1}(\gamma)$. For $v\in\widetilde{\gamma}$,
 set $\widetilde{v}:=\widetilde{\gamma}$.  Let $f|_A$ be the function obtained from $f$ by restricting its domain to $A$. 

For a positive integer  $n$, set $[n]:=\{1,\dots, n\}$. For a set  $X$ and element $y$, set $X+y:=X\cup \{y\}$ and  $X-y:=X\smallsetminus \{y\}$.

\section{Theorems related to bounded maximal degree}\label{maxd}

We begin this section with a short proof of \prettyref{thm:H-Sz}. This proof is
the basis for proofs of several stronger theorems, including the existence of an
efficient polynomial-time algorithm for equitably $k$-coloring { graphs when $k>\Delta$.
%We also consider conjectures and partial results related to bounded maximum degree. 
}

\subsection{Setup for equitable coloring in terms of maximum degree}

In this subsection, we organize a proof of \prettyref{thm:H-Sz} and many of its
generalizations. The notation is different, but the setup matches the original proof
until Case 2, where it takes a major short cut.

Suppose $G=(V,E)$ is a graph with $\Delta:=\Delta(G)$, $k$ is an integer with
 $k>\Delta$, and we want to equitably color $G$ with exactly $k$ colors. If $\left\vert G\right\vert =mk-p$,
where $1\leq p\leq k-1$, then set $G^{\prime}:=G+K_{p}$. Now $\left\vert G^{\prime}\right\vert $
is divisible by $k$, $\Delta(G^{\prime})\leq\max\{\Delta,k-2\}$, and the restriction
of any equitable $k$-coloring of $G^{\prime}$ to $G$ is an equitable $k$-coloring
of $G$. { So }%by induction
we may assume $n:=\left\vert G\right\vert =km$. Thus every
color class of an equitable $k$-coloring will have size $m$.

Argue by induction on $\|G\|$. If $G$ has no edges, we are done.  Else let $v\in V$ be non-isolated, and set $G':=G-E(v,V)$. By induction there  is an equitable $k$-coloring $f_{0}$ 
of $G'$. If $f_{0}$ is a proper coloring of $G$,   we are done. Else some
neighbors of $v$ are in $\widetilde{v}$. As $d(v)<k$, there is a class
$\widetilde{\beta}$ with no neighbors of $v$. Shift $v$ from   $\widetilde{v}$    %$\widetilde{x}$
to $\widetilde{\beta}$,  obtaining %leaving
  a new coloring $f$ with one \emph{small} class
$\widetilde{\alpha}_{0}:=\widetilde{v}-v$, one \emph{large} class $\widetilde{\beta}_{0}:=\widetilde{\beta}+v$,
and all other classes with size $m$. Such a coloring is  \emph{nearly equitable}.

To finish coloring  $G$, we will do a sequence of \emph{shift} operations that
improve our nearly equitable coloring, until finally it becomes equitable. First
define an auxiliary digraph $\mathcal{H=H}(f)$ on the colors of $f$ so that $\sigma\tau$
is a directed edge if and only if $\|y,\widetilde{\tau}\|=0$ for some vertex $y\in\widetilde{\sigma}$.\footnote{As the vertices of $\mathcal{H}$ are colors, not color classes, it is possible that
different colorings have the same $\mathcal{H}$.} In this case we say that $y$ \emph{witnesses} $\sigma\tau$. For edge $\sigma\tau\in\mathcal{H}$
with witness $w$, let \emph{$\sh(\widetilde{\sigma}\widetilde{\tau},w)$ }be the
operation that produces a new proper coloring by \emph{shifting} $w$ from $\widetilde{\sigma}$
to $\widetilde{\tau}$. We also write \emph{$\sh(\widetilde{\sigma}\widetilde{\tau},w)$
}for this coloring. When the particular witness is unimportant, we may simply write
$\sh(\widetilde{\sigma}\widetilde{\tau})$. We say that color $\sigma_{0}$ \emph{reaches}
color $\sigma_{k}$ and that $\sigma_{k}$ is \emph{reached by} $\sigma_{0}$ in
$\mathcal{H}$, if there is a directed $\sigma{}_{0}$,$\sigma_{k}$-path $\mathcal{P}=\sigma_{0}\dots\sigma_{k}\subseteq\mathcal{H}$.
In this case, put $\sh(\mathcal{P}):=\sh(\widetilde{\sigma}_{k-1}\widetilde{\sigma}_{k})\circ\dots\circ\sh(\widetilde{\sigma}_{0}\widetilde{\sigma}_{1})$.
This reduces the the size of $\widetilde{\sigma}_{0}$ by one, increases the size
of $\widetilde{\sigma}_{k}$ by one, and maintains all other class sizes. 

Let $\mathcal{A}$
be the set of colors that reach $\alpha_{0}$ in $\mathcal{H}$, and set $\mathcal{B}:=V(\mathcal{H})\smallsetminus\mathcal{A}$,
$a:=\left\vert \mathcal{A}\right\vert $, $b:=\left\vert \mathcal{B}\right\vert $,
{ $A:=\bigcup\{\widetilde{\alpha}:\alpha\in\mathcal{A}\}$, and $B:=V\smallsetminus A$.} So $k=a+b$ and $|B|=bm+1$.
%Now $\|v,\widetilde{\alpha}_{0}\|\geq1$ for all $v\in B$, so $bm+1\leq\|\widetilde{\alpha}_{0},B\|\leq\Delta(m-1)$.
%Thus $b<\Delta$ and
%\begin{equation}
%a\geq2.\label{eq:a>1}
%\end{equation}

Argue by a secondary induction on $b$, i.e., reverse induction on $a$.

\medskip{}

\noindent CASE 0: $\beta_{0}\in\mathcal{A}$; this includes the case $b=0$. Let
$\mathcal{P}$ be a $\beta_{0},\alpha_{0}$-path in $\mathcal{H}[\mathcal{A}]$.
Now $\sh(\mathcal{P})$ is an equitable $k$-coloring.

\medskip{}

Otherwise, $\beta_{0}\in\mathcal{B}$. A color $\alpha\in\mathcal{A}$ is \emph{terminal},
if every color $\gamma\in\mathcal{A}-\alpha$ reaches $\alpha_{0}$ in $\mathcal{H}-\alpha$.
Let $\mathcal{A}^{\prime}$ be the set of terminal colors, $a^{\prime}:=\left\vert \mathcal{A}^{\prime}\right\vert $,
and {$A':=\bigcup\{\widetilde{\alpha}:\alpha\in\mathcal{A'}\}$.} An edge $wz\in E(A,B)$ is \emph{solo}
if $f(w)\in\mathcal{A}^{\prime}$ and $\|z,\widetilde{w}\|=1$. In this case,  $w$
is a { \emph{solo root} of $z$%vertex
,}
 { and $z$ is a \emph{solo leaf} of $w$.} %HK: This was just a misprint from making the change. { Intuitively I don't think a vertex should be a solo leaf of itself, otherwise the following $\sh(\widetilde{y}\widetilde{\omega},y)$ makes less sense.}%neighbors of each other.}
For $z\in B$, let $s_{z}$ be the number of solo { roots (in $A'$) %neighbors
 of $z$}. Then
\begin{flalign}
\|z,A\| & \geq a+a'-s_{z}\geq a.\label{eq:ztoA}
\end{flalign}

\noindent CASE 1: { Some solo root %vertex
  $w\in\widetilde{\omega}$ 
  witnesses a solo edge $\omega\alpha$ of $\mathcal{H}[\mathcal{A}]$.} As $\omega\in\mathcal{A}'$,
some $\alpha,\alpha_{0}$-path $\mathcal{P}$  avoids $\omega$. Let $y$ be a solo { leaf %neighbor
 of $w$.} Set   
$f:=\sh(\widetilde{y}\widetilde{\omega},y)\circ\sh(\mathcal{P})\circ\sh(\widetilde{\omega}\widetilde{\alpha},w)$. Now $f|_{A+y}$
is equitable.   %and $f|_{G^*}$ is nearly equitable, where 
Set $G^{*}:=B-y$. By
\prettyref{eq:ztoA}, $\Delta(G^{*})<b$. By primary induction, $G^{*}$ has an equitable $b$-coloring $g$.  
 Thus $ f|_{A+y}\cup g$ is an equitable $k$-coloring of $G$. 

\medskip{}

Now assume no solo { root%vertex
} witnesses an edge of $\mathcal{H}[\mathcal{A}]$. 

\medskip{}

\noindent CASE 2: Some solo root $w\in\widetilde{\alpha}$ has two nonadjacent
solo { leaves%neighbors
 } $y,z\in B$; say $z\in\widetilde{\beta}$. 
 %Let $\cal{P}$ be an $\alpha,\alpha_{0}$-path in $\cal H$. Set $f:=\sh({\cal P})\cup \{(y,\alpha),(z,\alpha)\}-(w,\alpha)$.
   Set $B^{*}:=B-y$.
By \eqref{eq:ztoA}, $\Delta(B^{*})<b$. Using induction, equitably $b$-color
$B^{*}$. Now $\|w,B^{*}\|<b$, since $w$ witnesses no edge of $\mathcal{H}[\mathcal{A}]$
and $wy\in E$. Shift $w$ to some class   $\widetilde{\gamma}\subseteq B^{*}$, and
shift $y$ to $\widetilde{\alpha}-w$. This results in a new nearly equitable $k$-coloring
$f'$. Now $z$ witnesses that  $\beta\alpha\in\mathcal{H}(f')$, so $\beta\in\mathcal{A}(f')\supseteq\mathcal{A}$,
and we are done by secondary induction.

\subsection{Finishing the proof}

In this subsection we finish the proof of \prettyref{thm:H-Sz} and discuss how it
differs from the original proof. 
  Using Subsection 2.1, it suffices to show that one of Cases 0--2  holds. 
\begin{proof}[Proof of \prettyref{thm:H-Sz}]
 { Assume that Case 0 and Case 1 fail.} 
Then $b\geq1$ and no solo {root%vertex
} witnesses an edge of $\mathcal{H}[\mathcal{A}]$.
So for all solo {roots} $x$, 
\begin{equation}
\text{(a) }\|x,A\|\geq a-1~~\text{and~~(b) }\|x,B\|\leq d(x)-a+1\leq b.\label{eq:DegSolo}
\end{equation}
 Order $\mathcal{A}$ as $\alpha_{0},\dots,\alpha{}_{a-1}$ so that for all $i\in[a-1]$,
there is $j<i$ with $\alpha_{i}\alpha_{j}\in\mathcal{H}$; then $\alpha_{a-1}$
is terminal. 
Let $l$ be the least index such that $a_{j}$ is terminal whenever $l<j<a$.
%Let $\alpha_{l}$ be the last nonterminal color; by \prettyref{eq:a>1},
%$\alpha_{0}$ is a candidate for $l$. 
Now $a'\geq a-l-1$, and there is $t>l$ such
that all nontrivial $\alpha_{t},\alpha_{0}$-paths meet $\alpha_{l}$. 
As $\alpha_{t}$ has
no out-neighbors before $\alpha_{l}$, 
\begin{equation}
\text{(a) }\|w,A\|\geq l-1\geq a-a'-1\text{~~and~~(b) }\|w,B\|\leq b+a'\label{eq:nonsolo,A}
\end{equation}
 for all $w\in\widetilde{\alpha}_{t}$. Let $S$ be the set of solo {roots%vertices
}  in $\widetilde{\alpha}_{t}$,
and set $N:=N(S)\cap B$. Then $\|z,S\|\geq1$ for all $z\in N$ and $\|z,\widetilde{\alpha}_{t}\|\geq2$
for all $z\in B\smallsetminus N$. By \prettyref{eq:DegSolo}, $|N|\leq b|S|$, so
\begin{gather}
|N|+2|B\smallsetminus N|\leq\|S\cup(\widetilde{\alpha}\smallsetminus S),N\cup(B\smallsetminus N)\|\leq_{(\ref{eq:DegSolo},\ref{eq:nonsolo,A})} b|S|+(b+a')(m-|S|)\nonumber \\
2bm+2-b|S|\leq_{\prettyref{eq:DegSolo}} b|S|+bm-b|S|+a'(m-|S|)\nonumber \\
b(m-|S|)<a'(m-|S|)\nonumber \\
b<a'.\label{eq:b<a'}
\end{gather}

 Let $z\in B$. As $a+b-1\geq\Delta\geq d(z)\geq_{\prettyref{eq:ztoA}}a+a^{\prime}-s_{z}+\|z,B\|,$
\begin{gather}
s_{z}\geq a'-b+1+\|z,B\|.\label{eq:sig}
\end{gather}
Pick a maximum independent set $I\subseteq B$. Then $|I|+\sum_{z\in I}\|z,B\|\geq|B|=bm+1$.
Now
\begin{gather}
\sum_{z\in I}s_{z}\geq_{\prettyref{eq:sig}}\sum_{z\in I}(a'-b+\|z,B\|+1)>m(a'-b)+bm>_{\prettyref{eq:b<a'}}a'm=|A'|.\label{eq:C2}
\end{gather}
Thus some vertex in $A'$ has two (nonadjacent) solo {leaves%neighbors
} in $I$, and so Case
2 holds.
\end{proof}
This proof uses the failure of Case 1 to quickly derive $a'>b$. From this we easily
establish that Case 2 holds---some solo {root%vertex
} has two nonadjacent solo neighbors---and
quickly finish. In the original proof, a much more complicated argument used the
failure of Case 1 to derive $a'>2b$. From this another complicated argument proved
that Case 2 failed---no solo vertex has nonadjacent solo vertices. Finally, a nontrivial
effort established a contradiction. 

\subsection{A polynomial-time algorithm }

In 2007, Szemer\' edi mentioned to Kostochka that the original proof of \prettyref{thm:H-Sz}
does not provide a polynomial-time algorithm for finding the equitable coloring,
but that he had recently found one. Thinking about this led the first two authors
to a new proof of the theorem, another polynomial-time algorithm, and other new results
\cite{86,KKore}. Combining our ideas with those of Szemer\' edi and his student, Mydlarz, we  {developed} 
a simple, efficient algorithm for finding such a coloring.
\begin{thm}[Kierstead, Kostochka, Mydlarz and Szemer\' edi\cite{MR2676836}, 2010]
\label{thm:KKMS}Every graph on $n$ vertices with maximum degree less than $k$
can be equitably $k$-colored in $O(kn^{2})$ time using at most $2kn$ shifts.
\end{thm}

\begin{proof}
We will prove the bound on the number of shifts; the full result follows from standard
data-structure methods. For this we modify the proof of \prettyref{thm:H-Sz} in two ways. First, let $\mathcal{B}'\subseteq \mathcal{B}$ be the set of colors that can be reached from $\beta_{0}$,  { $B':=\bigcup\{\widetilde{\beta}:\beta\in\mathcal{B}'\}$,}
and $b':=|\mathcal{B}'|$. Now choosing $I\subseteq B'$ and using $\|z,B\|\ge\|z,B'\|+b-b'$ for $z\in B'\subseteq B$ in \prettyref{eq:C2}
 yields 
 \[%\begin{gather}
\sum_{z\in I}s_{z}\geq_{\prettyref{eq:sig}}\sum_{z\in I}(a'-b+\|z,B'\|+b-b'+1)>m(a'-b')+b'm>_{\prettyref{eq:b<a'}}a'm=|A'|,%\label{eq:C2}
\]%\end{gather}
so some solo root in $A'$ has two (nonadjacent) solo {leaves%neighbors
} in $I\subseteq B'$.
Thus we can obtain the first recoloring of $\mathcal{B}$
in Case 2 with at most $b'$ shifts. Secondly, we prove a stronger technical statement that will allow us to avoid the secondary induction on $b$. Instead we show by induction on $k$ that a nearly equitable $k$-coloring can be transformed to an equitable $k$-coloring using $2k-1$ shifts.
\begin{lem}
Let $G$ be a graph on $n=km$ vertices with a nearly equitable $k$-coloring $f$.
If 
\[\text{$d(x)<k$ for all $x\in A'\cup B,$}\tag{*}\]
 then $G$ can be equitably $k$-colored in
$2k-1$ shifts. 
\end{lem}

\begin{proof}
%Note that $a\ge2$: else $\mathcal A= \mathcal A'$, so $\Delta(G)<k$ and \prettyref{eq:a>1} holds. 
Argue as in the previous proof, but by induction on $k$. In Case 0,  $\sh(\mathcal{P})$
is an equitable $k$-coloring $G$ that uses $k-1$ shifts. In Case 1,
we use $a$ shifts to equitably $a$-color $G[A+y]$ and produce a nearly equitable
$b$-coloring of $G[B-y]$. Using induction, we %finish 
{equitably $b$-color $G[B-y]$}
in $2b-1$ shifts. So consider
Case 2. {As Case 0 fails,    \prettyref{eq:b<a'} implies} $1\le b<a'\le a$, so $a\ge 2$.    Now $\alpha_{0}$ is nonterminal, so $b<a'<a$. By integrality, $b+2\le a$.

Let $w$ be a solo root with nonadjacent solo leaves $y,z\in B'$; say $f(w)=\alpha$.
Then there is a $\beta_{0},f(y)$-path $\mathcal{P}\subseteq\mathcal{H}$. As Case
1 fails, some $x\in\widetilde{\alpha}-w$ is the first witness of an $\alpha,\alpha_{0}$-path
$\mathcal{Q}\subseteq\mathcal{H}$. Now $f':=\sh(\widetilde{y}(\widetilde{\alpha}-w),y)\circ\sh(\mathcal{Q},x)\circ\sh(\mathcal{P})$
is an equitable coloring of $G-w$ with a small class $\widetilde{\alpha}(f')=\widetilde{w}-w-x+y$.
We have used at most $a-1+b\leq2a-3$ shifts to obtain $f'$. Set $B^{*}:=B\cup\widetilde{\alpha}(f')$,
$G^{*}:=G[B^{*}]$, and $A^{*}=V\smallsetminus B^{*}$ . Every vertex $v\in B^{*}\smallsetminus\widetilde{\alpha}(f')$
satisfies $\|v,A^{*}\|\geq a-1$ (some vertices of $\widetilde{\alpha}(f')$
might not). So we can shift $w$ to some class $\widetilde{\beta}$ to obtain a nearly
equitable $(b+1)$-coloring $f^{*}$ of $B^{*}+w$ with big class $\widetilde{\beta}$. Now $z$  witnesses $\beta\alpha{\HK %\alpha does not change; the auxillary graph changes.
}\in{\cal H} 
(f')$, so {$\alpha$} is not terminal {in $f^*$}, and (*) holds {for $G^*$}. 
By induction, we can extend $f^{*}$ to an equitable $(b+1)$-coloring of $B^{*}$,
using {at most} $2b+1$ shifts. In total we have used {at most}  $(2a-3)+1+2b+1=2k-1$ shifts to equitably
$k$-color $G$.
\end{proof}

%Let $\sigma^{*}(k,n)$ be the number of shifts required to equitably $k$-color a graph $G$ on $n=km$ vertices with $\Delta(G)<k$.  For a nearly equitable $k$-coloring of $G$, let $\sigma(f)$ be the number of shifts required to convert $f$  to an equitable $k$-coloring, and set $\sigma(k)=\max\{\sigma(f):f\text{ nearly equitable}\}$.
%%$\sigma(a,b)$ be the number of shifts required to convert a nearly equitable $k$coloring $f$ with $a=a(f)$ and $b=b(f)$  to an equitable $k$-coloring.  
%The proof of Theorem~\ref{thm:H-Sz} yields:
%\[\sigma^{*}(k,n)\le\sigma^{*}(k,n-1)+1+\sigma(k),\]
%so it suffice to show that $\sigma(k)\le 2k-1$. 
%Our proof of \prettyref{thm:H-Sz} already provides an algorithm
%that uses at most $3nk$ shifts, since for each $i\in[n]$ it takes at most $3k$
%shifts to update an equitable $k$-coloring of $G_{i-1}$ 
% to a nearly equitable
%$k$-coloring of $G_{i}$, (one shift) and then to equitably $k$-color  $G_{i}$.
%To see this, argue by induction on $k$. In Case 1, we use $a$ shifts to decrease
%$k$ to $k-a$. Case 2 is more dangerous, since it involves the secondary induction
%twice, first when we update the nearly equitable $b$-coloring of $B_{0}:=B^{*}$
%and second when we update the nearly equitable $k$-coloring of $G_{i}$. But in
%Case 2, $a>a'>b$, so with care we can get the desired bound. We prefer a more structural
%approach that introduces two new ideas which somewhat improve the constant and may
%have other applications.

This completes the proof of \prettyref{thm:KKMS}. 
\end{proof}
\subsection{Equitable coloring in terms of maximum Ore-degree }\label{OreHS}

The Ore-degree of an edge $xy$ is the sum $\Theta(xy):=d(x)+d(y)$. The maximum
Ore-degree of a graph $G$ is $\Theta(G):=\max\{\Theta(xy):xy\in E\}$. Kostochka
and Yu \cite{KY0,KY4} conjectured that every graph $G$ with $\Theta(G)+1<2k$ is
equitably $k$-colorable; this would strengthen \prettyref{thm:H-Sz}. We proved  a
%an even 
stronger version. 

\begin{thm}[\cite{86}, 2008]
\label{thm:KK-Ore}If $\Theta{ (G)}<2k$ then $G$ has an equitable $k$-coloring. 
\end{thm}

\begin{proof}
Arguing by induction, we continue to use the same setup. In particular, we can assume
$|G|=km$ and there is a nearly equitable coloring $f$ of { $G'$}. For the latter,
{ let $G'=G-E(v,V)$, where $d(v)\le k$.} Also $\mathcal{H},\mathcal{A},\mathcal{A}',\mathcal{B}$,
etc.~are defined as before. These definitions depend on $f$, not any hypothesis
about maximum degree, { so \prettyref{eq:ztoA} and \prettyref{eq:DegSolo}(a) %, and \prettyref{eq:nonsolo,A}(a)
still hold.} %As in
{ By} \prettyref{eq:ztoA}, $\Theta(G[B])<2b$. We also need the subset
$\mathcal{A}''\subseteq\mathcal{A}'$ of colors that are separated from $\alpha_{0}$ { in $\mathcal{H}$}
by the last nonterminal color $\alpha_{l}$. Set $a'':=|\mathcal{A}''|$ and  { $A'':=\bigcup\{\widetilde{\alpha}:\alpha\in\mathcal{A}''\}$.} 
As in \prettyref{eq:nonsolo,A}{(a)}, for all $x\in A''$, 
\begin{equation}
{ \tag{{3*}}\|x,A\|\geq a-a''-1.\label{eq:4*}}%HK: Changed tag
\end{equation}
For $y\in B$ let { $s_{y}^{*}$} be the number of solo neighbors of $y$ in $A''$.
As in \prettyref{eq:ztoA}, 
\begin{equation}
{\|y,A\|\geq a+a''-s_{y}^{*}. \tag{{1*}}}\label{eq:2*}%HK: Changed tag
\end{equation}

Using the secondary induction, it is enough to find a nearly equitable $k$-coloring
$f'$ with $a(f')>a$. If Case 0 or Case 1 holds, we finish as before, so assume
both fail. The following lemma allows us to deal with Case 2 and to make progress,
even when all cases fail.
\begin{lem}
\label{lem:swap}Suppose Cases 0 and 1 fail and $wy$ is a solo edge {with $f(w):=\alpha\in\mathcal{A}''$}.
{ Then  $G'$ has a nearly equitable $k$-coloring $f'$} with $y\in \widetilde{\alpha}(f')-w$,
$f'(x)=f(x)$ for all $x\in A-w$, $\mathcal{A}\subseteq\mathcal{A}(f')$, and $\mathcal{A}''(f')\subseteq\mathcal{A}''$.
Moreover, if { $y'$ is another solo leaf of $w$ with $yy'\notin E$, then Case 0 holds for $f'$ or  $a(f')>a$.}
\end{lem}

\begin{proof}
Shift $y$ to $\widetilde{\alpha}-w$. 
{By  
induction}, { $B-y$} has an equitable
$b$-coloring $g$. By the failure of Case 1, \prettyref{eq:DegSolo}, and \prettyref{eq:ztoA},
$\|wy,A+y\|\geq2a$, so $\|w,B-y\|\leq2b-1$ . If $\|w,\widetilde{\beta}(g)\|=0$
for some $\beta\in\mathcal{B}$ then shift $w$ to $\widetilde{\beta}(g)$; else
$d(w)\geq a+b$, and $N(w)\cap\widetilde{\beta}(g)=\{z\}$ for some $\beta\in\mathcal{B}$
and $z\in\widetilde{\beta}(g)$. Then $d(z)<a+b$ and $\|z,A\|\geq a$, so $\|z,B\cup\widetilde{\alpha}-w\|\leq b-1$.
Thus $\|z,\widetilde{\gamma}(g)\|=0$ for some $\gamma\in\mathcal{B}+\alpha$ where
(i) $\widetilde{\gamma}(g)=\widetilde{\alpha}(g)=\widetilde{\alpha}+{y}-w$ or (ii)
$\gamma\in\mathcal{B}$. Shift $z$ to $\widetilde{\gamma}(g)$ and $w$ to $\widetilde{\beta}(g)$.
Let $f'$ be the final coloring obtained from $f$ and $g$. If (i)  holds then Case 0 holds
for $f'$. {Else if $wy'$ is solo and $yy'\notin E$, then $y'$ witnesses $f(y')\alpha\in \mathcal H(f')$, so $a(f')>a$.} %HK: Does it make sense now?
%ZX: This makes the last statement of this lemma a bit strange, because it seems that whether Case 0 holds or not has nothing to do with solo edge $wy'$ and $yy'\not\in E$.

Both $w$ and $y$ witness no edges of {$\mathcal{H}$} which are directed toward $\mathcal{A}$,
so $\mathcal{A\subseteq}A(f')$ and $E_{\mathcal{H}}(\mathcal{A}'',\mathcal{A}\smallsetminus\mathcal{A}'')=E_{\mathcal{H}(f')}(\mathcal{A}'',\mathcal{A}\smallsetminus\mathcal{A}'')$.
It follows that $\alpha_{l}$ separates $\mathcal{A}''$ from $\alpha_{0}$ in $\mathcal{H}(f')$.
But some terminal colors of $\mathcal{H}$ could be nonterminal colors in $\mathcal{H}(f')$.
Anyway, $\mathcal{A}''(f')\subseteq\mathcal{A}''$.
\end{proof}
We call the operation described in \prettyref{lem:swap} a $wy$-swap and write $f':=\sw(wy):= \sw_f(wy)$.
Assume Case 2 fails for $f$. In the remainder of the proof we show that some sequence
of swaps ends with {Case 0, 1 or 2.} Define a function $\mu:E(A'',B)\rightarrow\mathbb{Q}$
by $xz\mapsto\frac{b}{\|x,B\|}$, and for $z\in B$, let $\mu(z)=\sum_{x\in N(z)\cap A''}\mu(xz)$.
Now
\[
a''mb=|A''|b=\sum_{x\in A''}b=\sum_{xz\in E(A'',B)}\mu(xz)=\sum_{z\in B}\mu(z).
\]
By the Pigeonhole principle and the failure of Case 0, there is $y\in B$ with $\mu(y)<a''bm/|B|<a''$.
Pick $x\in A''\cap N(y)$ with $\|x,B\|$ maximum, and if $s_{y}^{*}\geq1$, let
$w\in A''$ be a solo neighbor of $y$ with $\|w,B\|$ maximum. So 
\begin{equation}
a''>\mu(y)\geq\frac{s_{y}^{*}b}{\|w,B\|}+\frac{2(a''-s_{y}^{*})b}{\|x,B\|}\text{, where }\|w,B\|:=1\text{ if }s_{y}^{*}=0.\label{eq:y*}
\end{equation}
Now (i) $\|w,B\|\geq b+1$ or (ii) $\|x,B\|\geq2b+1$: else $\mu(y)\geq s_{y}^{*}+(a''-s_{y}^{*})=a''$,
a contradiction. 
\begin{lem}
Vertex $y$ is incident to a solo edge from $A''$. 
\end{lem}

\begin{proof}
Suppose not. Then $s_{y}^{*}=0$, and by \prettyref{eq:y*}, $\|x,B\|\geq2b+1$.
This yields the contradiction 
\[
\Theta(xy)=d(x)+d(y)\geq_{\prettyref{eq:4*}\prettyref{eq:2*}}(a-a''-1+2b+1)+(a+a'')=2a+2b=2k.\qedhere
\]
\end{proof}
First suppose $a''>b$. Set $f':=\sw(wy)$. We claim that  one of the cases holds
for $f'$ or (iii)~$\|A''(f'),B(f')\|<\|A'',B\|$. Suppose (iii) fails. Then $\|w,B\|+\|y,A''\|\leq\|w,A''\|+\|y,B\|$.
Counting $wy$ twice yields, 
\begin{gather*}
2k\geq\|wy,(A\smallsetminus A'')\cup A''\cup B\|\geq2(a-a'')+\|w,B\|+\|y,A''\|+\|w,A''\|+\|y,B\|\\
a''+b\geq\|w,B\|+\|y,A''\|.
\end{gather*}
Now $\|y,A''\|\geq a''$, so (iv)~$\|w,B\|\leq b$. Thus (v)~$s_{y}^{*}\geq2a''-\|y,A''\|\geq a''-b+\|w,B\|$.
So
\[
\mu(y)\geq_{\prettyref{eq:y*}}\frac{s_{y}^{*}b}{\|w,B\|}\geq_{\text{(v)}}\frac{a''-b+\|w,B\|}{\|w,B\|}b\geq(a''-b)\frac{b}{\|w,B\|}+b\geq_{\text{(iv)}}a'',
\]
a contradiction. So each unsuccessful swap decreases $\|A'',B\|$. Eventually, one
of the cases holds or $a''\leq b$. 

Now suppose $a''\leq b$. Then using $\prettyref{eq:y*}$, $\|w,B\|>s_{y}^{*}$.
Set $f':=\sw(wy)$. We claim that $a<a(f')$ or one of the cases holds for $f'$
or $\|B\|<\|B(f')\|$. Suppose not; then $\|B\|\geq\|B(f')\|$. Now $\|y,B\|\geq\|w,B\|-1\geq s_{y}^{*}$.
Now we have the contradiction:
\begin{align*}
\Theta(G)\geq\begin{cases}
\Theta(wy)\geq_{{(\ref{eq:DegSolo},\ref{eq:ztoA})}}(a-1+\|w,B\|)+(a+\|y,B\|)\geq2k & \text{if (i) {holds,}}\\
\Theta(xy)\geq_{{(\ref{eq:4*},\ref{eq:DegSolo})}}(a-a''-1+\|x,B\|)+(a+a''-s_{y}^{*}+\|y,B\|)\geq2k & \text{if (ii) {holds.}}
\end{cases}
\end{align*}
By \prettyref{lem:swap}, $a''(f')\leq a''\leq b=b(f')$. So each unsuccessful swap
increases $\|B\|$.
\end{proof}
Because of the repeated use of \prettyref{lem:swap}, this proof only yields an exponential-time
algorithm.
\begin{conjecture}
There is a polynomial-time algorithm for finding an equitable $k$-coloring of every
graph with maximum Ore-degree $\Theta(G)<2k$. 
\end{conjecture}

\subsection{Equitable coloring of directed graphs}

Next we consider the problem of extending the Hajnal-Szemerdi Theorem to directed
graphs. It is not immediately clear what this should mean. Let $G$ be a \emph{simple}
directed graph---a loopless digraph such that 
 for every ordered pair $(x,y)$ of vertices there is at most one edge from $x$ to $y$.
%there are at most two edges $xy$ and $yx$ between any two vertices $x$ and $y$. 
Let $d^{-}(v)$ and $d^{+}(v)$
denote the in- and out-degree of $G$. The \emph{maximum} and \emph{minimum} degree
of $G$ are $\Delta(G):=\max\{d^{-}(v)+d^{+}(v)\}$ and $\delta(G):=\min\{d^{-}(v)+d^{+}(v)\}$.
A directed graph is \emph{acyclic} if it contains no directed cycle. In the context
of directed graphs, an \emph{acyclic coloring }is a coloring whose classes are acyclic
(instead of independent). 
\begin{thm}[Czygrinow, DeBiasio, Kierstead  and  Molla \cite{113}, 2015]
\label{thm:dH-Sz}Every directed graph $G$ with maximum degree $\Delta(G)<2k$
has an equitable acyclic $k$-coloring.
\end{thm}

The theorem implies the Hajnal-Szemer\' edi Theorem. Suppose $G=(V,E)$ is an undirected
graph with $\Delta(G)<k$. Let $\vec{\ensuremath{G}}$ be the directed graph obtained
from $G$ by replacing every edge $uv\in E$ by directed edges $uv$ and $vu$, and
observe ({*}) $vuv$ is a cycle in $\vec{G}$. Now $\Delta(\vec{G})\leq2k-2<2k-1$,
so $\vec{G}$ has an equitable acyclic $k$-coloring $f$ by \prettyref{thm:dH-Sz}.
Using ({*}), $f$ is an equitable $k$-coloring of $G$.

Here is another consequence of \prettyref{thm:dH-Sz}. A simple directed graph $G$
is complete if it contains every possible edge. If $G$ is complete and $|G|=m$,
we write $G=\vec{K}_{m}$. Set $\overline{G}:=\vec{K}_{|G|}-E(G)$. An \emph{equitable
$k$-tiling} of a directed graph $G$ on $km$ vertices by transitive tournaments
is a $k$-partition of $V(G)$ such that each part induces a transitive tournament. 
\begin{cor}
\label{cor:TransTour}Every directed graph $H$ on $n=km$ vertices with $\delta(H)\geq2n(1-\frac{1}{m})-1$
has an equitable $k$-tiling by transitive tournaments.
\end{cor}

\begin{proof}
Set $G:=\overline{H}$. Then $\Delta(G)=2n-\delta(H)-1\leq2k-1<2k$. By \prettyref{thm:dH-Sz},
$G$ has an equitable acyclic coloring $f$, so $G[\widetilde{\gamma}]$ is acyclic { for every class $\widetilde{\gamma}$.} 
Thus $G[\widetilde{\gamma}]$ is a spanning subgraph of a transitive tournament $T$.
So $\overline{T}$ is a transitive tournament spanning $\overline{G}[\widetilde{\gamma}]=H[\widetilde{\gamma}]$.
\end{proof}
The proof of Theorem~\ref{thm:dH-Sz} is very similar to that of Theorem~\ref{thm:H-Sz},
and it even provides a poynomial-time algorithm. One significant difference is in
the definition of the auxiliary graph $\mathcal{H}$. Now $\gamma\alpha\in\mathcal{H}$
if and only if for some $x\in\widetilde{\gamma}$, there is no cycle in $\widetilde{\alpha}+x$.
In particular, if $\gamma\alpha\notin\mathcal{H}$ then for all $x\in\widetilde{\gamma}$
there are edges from $x$ to $\widetilde{\alpha}$ and from $\widetilde{\alpha}$
to $x$.

One might ask why, in the statement of Corollary~\ref{cor:TransTour}, attention
is restricted to transitive tournaments. Of course one answer is that they are easy
to handle. The special (and exceptional) case $m=n/k=3$ provides some direction.
There are two tournaments on three vertices: one $\vec{T}$ is transitive, the other
$\vec{C}$ is cyclic. The next theorem gives a degree condition for tiling $\vec{G}$
with various mixtures of transitive and cyclic tournaments.
\begin{thm}[Czygrinow, Kierstead  and  Molla \cite{107}, 2014]
\label{thm:T-packing}Suppose $\vec{G}$ is a directed graph on $n=3k$ vertices
with $\delta(\vec{G})\geq\frac{4}{3}n-1$. For all $t\in[k]$, $\vec{G}$ has a $k$-tiling
by $t$ transitive and $k-t$ cyclic tournaments.
\end{thm}

Notice that the theorem does not assert the existence of a tiling entirely by cyclic
tournaments since $t\geq1$ is assumed. This is necessary. Indeed, in 2000 Wang~\cite{Wangdir},
observed that even $\delta(\vec{G})\geq\frac{3|\vec{G}|-5}{2}$ is not enough to
force a tiling with cyclic tournaments. Let $3$ divide $2s+1$, and { set $G:=\vec{K}_{s}+\vec{K}_{s+1}+\vec{K}_{s,s+1}$, where $\vec{K}_{s,s+1}$ denotes the complete bipartite graph  $K_{s,s+1}$ with all edges oriented from the small part to the big part.} 
Then $\vec{G}$ satisfies the degree condition. As one of $s$ and $s+1$ is not
divisible by $3$, any $3$-tiling has a $3$-tournament that meets both $\vec{K}_{s}$
and $\vec{K}_{s+1}$. Any such tournament must be transitive because its two crossing
edges are oriented in the same direction. Wang also proved that this example is tight. 

A multigraph $G$ is \emph{standard} if it is loopless and has multiplicity at most
$2$. An $s$\emph{-sided triangle} is a standard multigraph on three vertices that
contains $K_{3}$ and has $s$ edges. To prove Theorem~\ref{thm:T-packing}, first 
 { obtain} a standard multigraph $G$ by removing the directions from all edges of
$\vec{G}$. Then the minimum degree of $G$ as a multigraph is the same as the minimum
degree of $\vec{G}$ as a simple directed graph. Any $5$-sided triangle of $G$
spans both a cyclic and a transitive tournament in $\vec{G}$, and a $4$-sided triangle
of $G$ spans a transitive tournament in $\vec{G}$. Thus Theorem~\ref{thm:T-packing}
follows from:
\begin{thm}[Czygrinow,  Kierstead and Molla~\cite{107}, 2014]
Suppose $G$ is a standard multigraph on $n=3k$ vertices with $\delta(G)\geq\frac{4}{3}n-1$.
Then $G$ has a $k$-tiling by triangles such that all but one of the triangles are
$5$-sided, and the remaining triangle is $4$-sided. 
\end{thm}

Call a directed graph $H$ universal if $T\subseteq H$ for every tournament $T$
with $|T|=|H|$.
{ \begin{conjecture}[Czygrinow, DeBiasio, Kierstead and Molla~\cite{113}, 2015]
\label{conj:univ} Suppose $m\ne3$ and $G$ is a directed graph on $n=km$
vertices with $\delta(\vec{G})\geq2(1-\frac{1}{m})n-1$. Then $\vec{G}$ has a $k$-tiling
by universal directed graphs on $m$ vertices. 
\end{conjecture}} 

\begin{thm}[Czygrinow, DeBiasio, Molla and Treglown \cite{CDMT18}, 2018]
{ For all $m$, the preceding conjecture holds for sufficiently large $n$.} 
\end{thm}

%\section{}

\section{Chen-Lih-Wu Conjecture}

Brooks' Theorem states that every graph $G$ with $\Delta(G)\leq k$ is $k$-colorable
except when $K_{k+1}\subseteq G$ or $k=2$ and $G$ contains an odd cycle. Could
we also require that the $k$-coloring be equitable? When $k$ is odd, $G=K_{k,k}$
is a $2$-colorable graph with $\Delta(G)=k$ and $|G|=2k$. Thus any equitable $k$-coloring
of $G$ consists of $k$ classes of size $2$, but each color class is contained
in one part of $G$, so $G$ has at most $2\lfloor\frac{k}{2}\rfloor=k-1$ such classes. So, $K_{k,k}$ is not equitably $k$-colorable when $k$ is odd.
The following conjecture, asserting that this is the only counterexample, is the
most important open question on equitable coloring.
\begin{conjecture}[Chen, Lih and Wu \cite{CLW}, 1994]
\label{conj:CLW}If $G$ is a $k$-colorable graph with $\Delta(G)\le k$, then
either $G$ is equitably $k$-colorable or $k$ is odd and $G$ contains a $K_{k,k}$.
\end{conjecture}

By \prettyref{thm:H-Sz}, we only need to consider the case $k=\Delta(G)$. 
Lih and Wu~\cite{MR1391262} proved the conjecture for bipartite graphs.
For general
graphs, the conjecture is only known to hold for  special cases of $k$. The
case $k\leq2$ is trivial. Chen, Lih and Wu \cite{CLW} themselves proved
the cases $k=3$ and $k\geq |G|/2$. { The first two authors} proved the following.
%Theorem.
\begin{thm}
\label{thm:PreCLW}Suppose $G$ is a $k$-colorable graph with $\Delta(G)=k$ such
that $K_{k,k}\nsubseteq G$ if $k$ is odd. Then $G$ has an equitable coloring if~$k\leq4$ \emph{(\cite{KK11}, 2012)} or~$k\geq|G|/4$ \emph{({\cite{109},} 2014)}.
\end{thm}

{ The arguments used in~\cite{KK11} and~{ \cite{109}} are based on the approach described in Section 2,} but every step is more difficult. For example, the existence of a nearly equitable coloring in the proof of Hajnal-Szemer\' edi Theorem is obvious, but to show this existence in the proof of Theorem~\ref{thm:PreCLW} we used a sharpening of Brooks' Theorem due to Kostochka and Nakprasit~\cite{KN05}.  They showed that if $\Delta(G),\chi(G)\le k$ then for every  $k$-coloring $f$ of $G-v$ there is a $k$-coloring $f'$ of $G$ such that for all colors $\gamma$, except one,   $|\widetilde{\gamma}(f)|=|\widetilde{\gamma}(f')|$. Another difficulty is that while { using induction in the proof of Hajnal-Szemer\' edi Theorem,  we automatically have an equitable $b$-coloring of $G[B-y]$, but in the proof of Theorem~\ref{thm:PreCLW} we need to check that  when $b$ is odd,  $G[B-y]$ contains no complete    $K_{b,b}$. We will have more to say about this in the next subsection.}

It is notable that the state of knowledge (\prettyref{thm:PreCLW}) concerning Conjecture
\ref{conj:CLW} is analogous to that (\prettyref{thm:PreHS}) concerning Erd\H{o}s'
conjecture before the proof of \prettyref{thm:H-Sz}.

For sparse graphs, more is known. Yap and Zhang \cite{YZ1} proved that the conjecture
holds for all outerplanar graphs.
Yap and Zhang \cite{YZ2} also proved that the conjecture
holds for planar graphs when $k\geq13$, Nakprasit \cite{KitThesis,Na12}
proved it for planar graphs when $9\leq k\leq12$, and
Kostochka, Lin and Xiang \cite{KLX} { settled the case $k=8$.}
%{\HK gave a unified proof for planar graphs with maximum degree at least $8$.???}
Thus, for planar graphs Conjecture~\ref{conj:CLW} is not verified only for graphs with maximum degree $5,6$
and $7$. 
For these degrees,
 Conjecture~\ref{conj:CLW}  was proved for planar graphs with extra restrictions, mainly with restrictions on cycle structure. 
%In particular, these rese arch considered planar graphs without certain short cycles.
In 2008, Zhu and Bu~\cite{ZhuBu} proved that the conjecture holds for $C_3$-free planar graphs with maximum degree $\Delta\geq 8$. It also holds for $C_4,C_5$-free planar graphs with maximum degree $\Delta\geq 7$. In 2009, Li and Bu~\cite{LiBu} proved that the conjecture holds for $C_4,C_6$-free planar graph with maximum degree $\Delta\geq 6$. In 2012, Nakprasit and Nakprasit~\cite{NaNa} proved that the conjecture holds for $C_3$-free planar graphs with maximum degree $\Delta\geq 6$, $C_4$-free planar graphs with maximum degree $\Delta\geq 7$, and planar graphs with maximum degree $
\Delta\geq 5$ and girth at least $6$.

Zhang~\cite{1planar} proved the conjecture for the wider class of $1$-planar graphs, that is, graphs that admit a drawing in the plane such that each edge crosses at most one other edge
 but with the stronger restriction on $k$: for  $k\geq 17$. 
The restriction $k\geq 17$ was relaxed to $k\geq 13$ for $1$-planar graphs by Cranston and Mahmoud~\cite{CM24}. In fact, Zhang %(but not Cranston and Mahmoud)
proved his bound $17$ 
 for the more general class $\mathcal{Z}$ of graphs $G$ such that $\|G[A]\|\leq 4(|A|-2)$
for every $A\subseteq V(G)$ with $|A|\geq 3$.

\subsection{An equivalent conjecture}
\begin{figure}[ptbh]
\begin{center}
\includegraphics[scale=0.45]{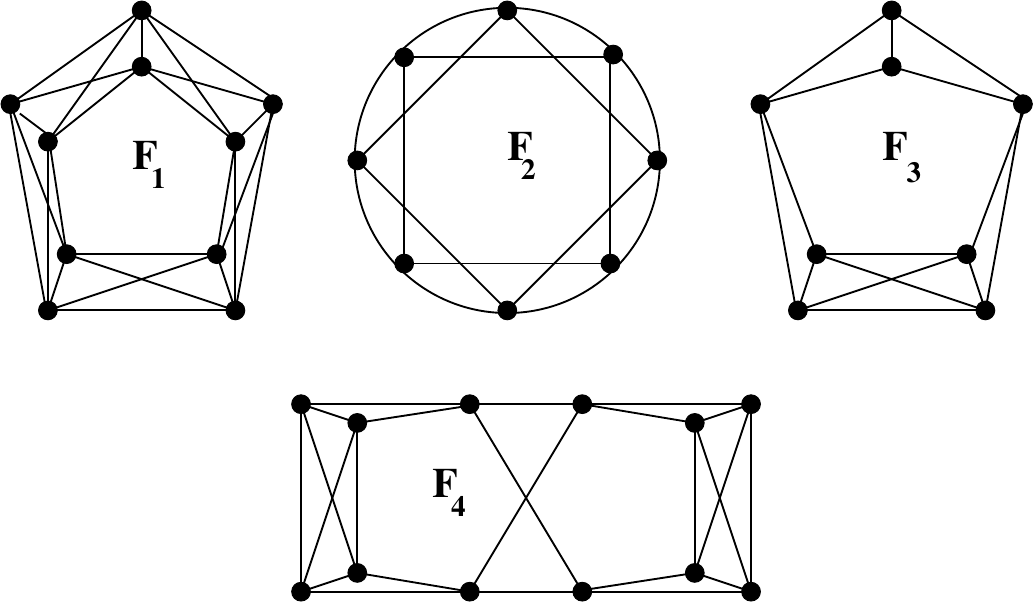}
\end{center}
\caption{One $5$-equitable and three $4$-equitable basic graphs.}%
\label{f2}%
\end{figure}

\begin{figure}[ptbh]
\begin{center}
\includegraphics[scale=0.37]{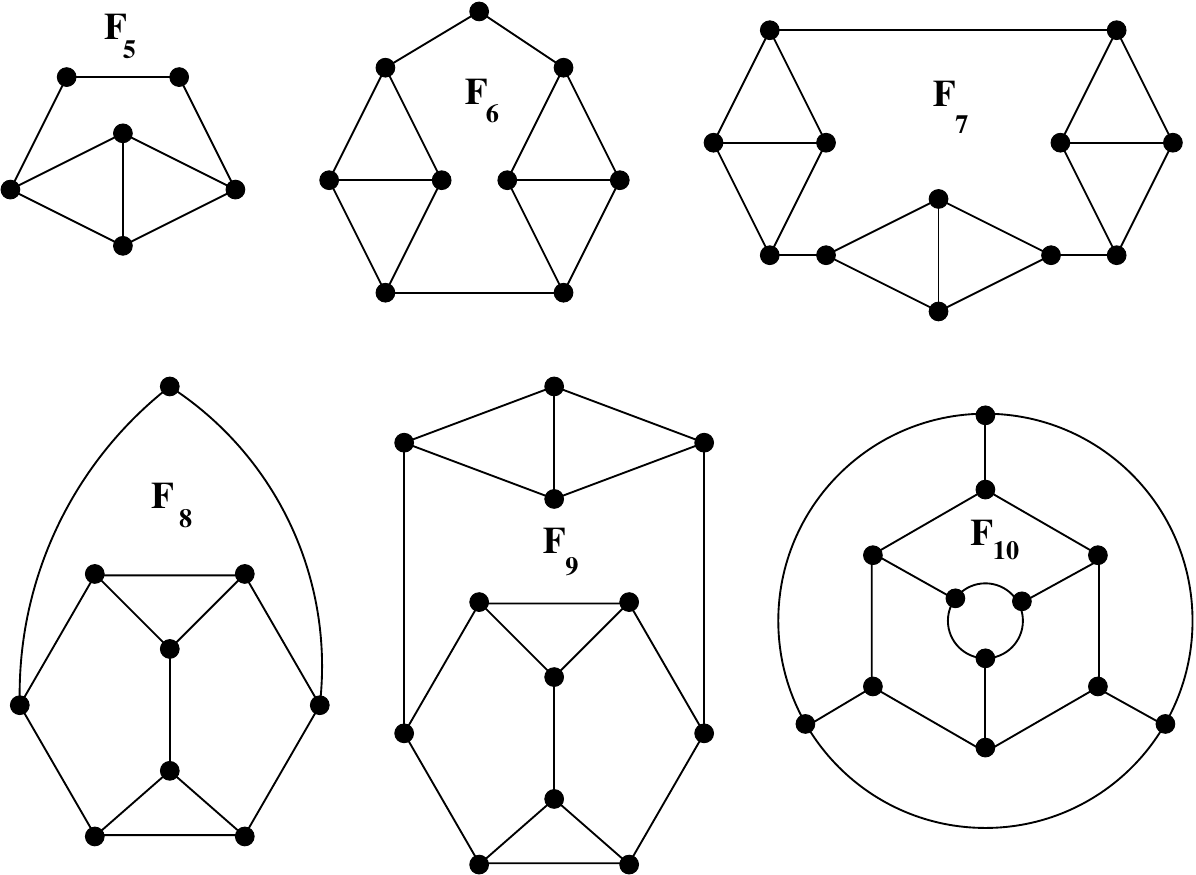}
\end{center}
\caption{Six $3$-equitable basic graphs.}%
\label{f3}%
\end{figure}

Note that for an odd $k$, Conjecture~\ref{conj:CLW} does not describe all
graphs $G$ with $\Delta(G)\leq k$ that are
not equitably $k$-colorable. For example,
for each odd $k\geq3$, the graph consisting of two disjoint copies of $K_{k,k}$
has an equitable $k$-coloring, but the graph consisting of a copy of
$K_{k,k}$ and a copy of $K_{k}$ does not. This construction can be generalized. Say
that a graph $H$ is ${k}$\emph{-equitable }if $\left\vert H\right\vert $
is divisible by $k$, $H$ is $k$-colorable and every $k$-coloring of $H$ is
equitable. For an odd $k$, if $G$ contains $K_{k,k}$ and $G-K_{k,k}$ is $k$-equitable, then $G$
does not have an equitable $k$-coloring. This motivates the study of equitable graphs,
i.e., graphs that are $k$-equitable for some $k$.
It was proved in~\cite{KKclw} 
that  there is a good
description of the family of all $k$-equitable graphs; they can all
be built from simple examples in a simple way.
Note first that
if a spanning subgraph of a graph $G$
is the disjoint union of $k$-equitable graphs, then $G$
is $k$-equitable. Clearly, $K_k$ is $k$-equitable. Figures~\ref{f2}
and~\ref{f3} show $10$ other equitable $2$-connected graphs.

Together with $K_{k}$, the
$k$-equitable graphs in Figures~\ref{f2}
and~\ref{f3}  are the { $k$\emph{-basic}} graphs.
A
{ $k$\emph{-decomposition}} of $G$ is a partition on $V(G)$ into subsets
$V_{1},\dots,V_{t}$ such that each $G[V_{i}]$ is { $k$-basic}.

\begin{thm}[\cite{KKclw}, 2010]
\label{First}
Let $k\geq3$. Let $G$ be an { $k$-colorable} graph with
$\Delta(G)= k$ and $|V(G)|$ divisible by $k$.
Then the following are equivalent:
\begin{enumerate}[label=(\alph*),nosep]
\item  $G$ has
a { $k$-decomposition};
\item  $G$ is { $k$-equitable};
\item  $G$ has an equitable $k$-coloring but does not have a nearly
equitable $k$-coloring.
\end{enumerate}
\end{thm}

 Theorem~\ref{First} is used in the proof of Theorem~\ref{thm:PreCLW} to get $b$-equitable colorings of $B-y$. It also  inspired the following conjecture which (if true) would describe all
graphs $G$ with $\Delta(G)\leq k$ that are
not equitably $k$-colorable.

\begin{conjecture}[\cite{KKclw}, 2010]
\label{strong}
Suppose that $k\geq3$ and $G$ is a {$k$-colorable} graph with
$\Delta(G)= k$. Then $G$ has no equitable $k$-coloring if and only if $k$
is odd and there exists $H\subseteq G$ such that { $H=K_{k,k}$} and $G-H$ has
a { $k$-decomposition}.
\end{conjecture}

Theorem~\ref{First} yields the following implications.

\begin{cor}\label{kkc1}
  For all positive integers $k$ and $n>k$,
Conjecture~\ref{strong}
holds for all $k$-colorable graphs $G$ with {  $\Delta(G)\leq k$ and
at most $n$ vertices} if and only if
Conjecture~\ref{conj:CLW} holds for all such graphs.
\end{cor}
%\begin{conjecture}[] \end{conjecture}

\begin{cor}\label{kkc2}
Every $k$-equitable $k$-colorable graph $G$ with $\Delta(G)\leq k$  has a unique $k$-decomposition.
\end{cor}

\begin{cor}\label{kkc3}
  There exists a polynomial time algorithm for deciding whether a
  graph $G$ with $\Delta(G)\leq k$
 is $k$-equitable.
\end{cor}

\subsection{An Ore version of the Chen-Lih-Wu Conjecture}

{Recall that the Ore-degree of  a graph $G$ is $\Theta(G):=\max\{d(x)+d(y):xy\in E\}$.}
 The analog of the trivial bound $\chi(G)\leq \Delta(G)+1$ in terms of the Ore-degree is {
 the  claim that if % $\chi(G)\leq \left\lfloor\frac{\Theta(G)}{2}\right\rfloor+1$. 
  $2k>\Theta(G)$, then $G$ is $k$-colorable.}
 This claim is true and is easy to prove. The bound is achieved for complete graphs and odd cycles. Moreover, there are infinitely many graphs $G$ with 
$\Theta(G)=7$ and $\chi(G)=4$. Also the $5$-Ore graph with $9$ vertices has Ore-degree $9$ and chromatic number $5$. However, the following analog of Brooks' Theorem for the Ore-degree holds.

\begin{thm}[\cite{KK09}, 2009]
\label{Ore2}
If $k\geq 7$ and $G$ is a connected graph with $\Theta(G)\leq 2k+1$ distinct from $K_{k+1}$, then $\chi(G)\leq k$.
\end{thm}

We think that the claim of Theorem~\ref{Ore2} holds also for $k=6$ but cannot prove this. 

{Theorem ~\ref{thm:KK-Ore} in Section~\ref{OreHS} proves an analog of Theorem ~\ref{thm:H-Sz} in terms of the Ore-degree.
So,} in view of Theorems~\ref{Ore2} and~\ref{thm:KK-Ore}, the following 
{analog of Conjecture~\ref{conj:CLW}} was posed in~\cite{KKore} and~\cite{KK09}.

\begin{conjecture}[\cite{KKore,KK09}, 2009]
\label{Ore-s2}
Let $G$ be a connected graph with $\Theta(G)\leq 2k$ and $k\geq 3$. If $G$ is distinct from $K_{k+1}$ and from $K_{m,2k-m}$
for odd $m$, then $G$ has an equitable $k$-coloring.
\end{conjecture}

This conjecture was proved for $k=3$. In~\cite{KK09}, the following strengthening of it was conjectured. 

\begin{conjecture}[\cite{KK09,KKY}, 2009]
\label{Ore-s3}
Suppose that $k\geq 3$. A $k$-colorable, $n$-vertex graph $G$ with 
$\Theta(G)\leq 2k$ has no equitable 
$k$-coloring, if and only if $n$ is divisible by $k$, and there exists 
$W\subset V(G)$ such that $G[W]=K_{m,2k-m}$ for some odd $m$
and $G-W$ is $k$-decomposable.
\end{conjecture}

It was proved in~\cite{KK09} 
that Conjectures~\ref{Ore-s2} and~\ref{Ore-s3} are  equivalent even for graphs with restricted number of vertices and restricted values of Ore-degree. This implies that Conjecture~\ref{Ore-s3} holds for $k=3$. This stronger claim 
may help
 to prove Conjecture~\ref{Ore-s2} by induction, because we would need to prove the formally weaker conjecture but use the fact that the stronger conjecture holds for smaller graphs.

%\subsubsection{The conjecture, examples and partial results}

%\subsubsection{An equivalent conjectured}

%\subsection{Stronger and stronger conjectures}

\iffalse
\section{Strongly equitable list-coloring}

Small and large degree bounds. Stronger Conjecture

\subsection{Examples}

\subsection{Degree bounds}

Check that old bounds work for SE

\fi

\section{Equitable list coloring}

Returning to the initial idealized example of this paper about equitable distribution of jobs among workers, % we can consider 
 %a somewhat 
we consider the more 
 complicated situation  when not all workers have the same skills. Maybe each job $J$ has a list $L(J)$ of workers who can perform this job. Then it might be that some workers are able to perform only few jobs. In this case, it is not realistic to insist on almost the same number of jobs being performed by each worker. It is more natural to ask for an upper bound on the number of jobs performed by each worker.
Among the models for this more complicated situation there are different kinds of {\em equitable list coloring}.

A \emph{list assignment} $L$ for a graph $G$ assigns to each vertex $v$ a set $L(v)$ of  colors.
An \emph{$L$-coloring} $f$ is a proper coloring of $G$ such that $f(v)\in L(v)$ for each $v\in V(G)$.
If $|L(v)|=k$ for each $v\in V(G)$, then $L$ is a \emph{$k$-list assignment}. 

{Kostochka, Pelsmajer and West~\cite{KPW} proposed a list version of equitable coloring as follows.}
For a $k$-list
assignment $L$, a graph $G$ is \emph{equitably $L$-colorable} if it has an $L$-coloring
such that every color class has %size
at most $\lceil|G|/k\rceil$ {vertices}. {Now} $G$ is \emph{equitably
$k$-choosable} if %$G$
{it} is  equitably $L$-colorable for every $k$-list assignment
$L$.

In contrast to ordinary equitable coloring, {in this setting} it is possible that the sizes of two color classes differ {significantly}, %by more than $1$, 
as {it could be that} some color  {appears only} %show up
 in the lists of  few vertices.

The following analog of Theorem~\ref{thm:H-Sz} was  conjectured in~\cite{KPW} .
\begin{conjecture}[Kostochka, Pelsmajer and West, Conjecture 1.1 \cite{KPW}, 2003]\label{conj:list-HS}
    Every graph $G$ is equitably $k$-choosable for each $k> \Delta(G)$.
\end{conjecture}

 %They~\cite{KPW} confirmed Conjecture~\ref{conj:list-HS} when $k\ge\max\{\Delta(G)+1,|G|/2\}$.
 %HK: I do not like starting a paragraph with a pronoun. Maybe: In ~\cite{KPW}  Conjecture~\ref{conj:list-HS} was confirmed for $k\ge\max\{\Delta(G)+1,|G|/2\}$.
%
%and Conjecture~\ref{conj:list-CLW} when $k\ge \max\{\Delta(G), |G|/2\}$.
In~\cite{KPW}  Conjecture~\ref{conj:list-HS} was confirmed for $k\ge\max\{\Delta(G)+1,|G|/2\}$. 
Pelsmajer~\cite{MP} and Wang and Lih~\cite{WL} independently proved Conjecture~\ref{conj:list-HS} for $k\le 3$.

In 2013~\cite{KK13} the conjecture was proved for a larger class of graphs; several relaxations of the conjecture were also proved.

\begin{thm}[\cite{KK13}, 2013]\label{thm:KK13}
    Let $G$ be a graph with $\Delta(G)=\Delta$.
    \begin{enumerate}[label=(\alph*),nosep]
        \item If $\Delta\le 7$ and $k\ge \Delta+1$, then $G$ is equitably $k$-choosable.
        \item If $|G|\ge \Delta^3$ and $k\ge \Delta+2$, then $G$ is equitably $k$-choosable.
        \item If $\omega(G)\le k$ and $|G|\ge 3(k+1)k^8$, then $G$ is equitably $(k+1)$-choosable.
    \end{enumerate}
\end{thm}

For each of the conditions in Theorem~\ref{thm:KK13}, there is a polynomial-time (in terms of $n$) algorithm, that finds an equitable $L$-coloring of $G$ for every graph $G$ satisfying the condition, and every $k$-list assignment $L$ on $G$. The method used to prove Theorem~\ref{thm:KK13} is based on ideas of the proof of the  Hajnal-Szemer\' edi Theorem and its variations, such as
Theorem~\ref{thm:PreCLW}. But several new tricks were applied to handle {distinct} %different
lists. In particular, the definition of the auxiliary digraph $\mathcal{H}$ needed to be modified {so that $\sigma\tau$ is a directed edge of $\cal H$ if and only if there is $w\in \widetilde{\sigma}$ with $\tau\in L(w)$ and $\|w,\widetilde{\tau}\|=0$.} Still, these tricks were insufficient to  prove the conjecture in full.

%HK: As previously stated, each graph could have its own algorithm, making the theorem trivial. 
%For each graph $G$ that satisfies one of the conditions in Theorem~\ref{thm:KK13}, there is a polynomial time algorithm, in terms of $n$, that finds an equitable $L$-coloring for every $k$-list assignment $L$ on $G$.

An Ore-type version of Conjecture~\ref{conj:list-HS} 
(which also can be viewed as the list version of Theorem~\ref{thm:KK-Ore}) 
was proposed in~\cite{KKY}:

\begin{conjecture}[Conjecture 6.3 \cite{KKY}, 2009]\label{conj:Orelist}
    Every graph $G$ is equitably %{\HK Coordinate} 
    $(1+0.5\Theta(G))$-choosable.
\end{conjecture}

In support of this conjecture, it was proved in~\cite{KKY} that every graph $G$ with $\Theta(G)\leq 6$ is equitably $4$-choosable. 
{
A Chen-Lih-Wu version of Conjecture~\ref{conj:list-HS}  {(in other words, the list version of the Chen-Lih-Wu Conjecture)}
was also  proposed in~\cite{KPW}} in the following form.
%The CLW version of Conjecture~\ref{conj:list-HS} has also been proposed in~\cite{KPW} in the following form.

\begin{conjecture}[Conjecture 1.2, Kostochka, Pelsmajer and West \cite{KPW}, 2003]\label{conj:list-CLW}
    If $G$ is an $k$-colorable graph with $3\le\Delta(G)\le k$, then either $G$ is equitably $k$-choosable or $k$ is odd and $G$ contains a $K_{k,k}$.
\end{conjecture}

Conjecture~\ref{conj:list-CLW} was confirmed for $k\ge\max\{\Delta(G),|G|/2\}$, connected interval graph and $2$-degenerate graphs $G$ with $k\ge\max\{\Delta(G),5\}$ in~\cite{KPW}.  Zhang and Wu~\cite{ZhW} proved the Conjecture  for 
series-parallel graphs.
Dong and Zhang~\cite{DZ18} proved that if $G$ is a graph with $\mad(G)<3$ then $G$ is equitably
$k$-choosable for $k\ge\max\{\Delta(G),4\}$, and if $\mad(G)<\frac{12}{5}$ then
$G$ is equitably $k$-choosable for $k\ge\max\{\Delta(G),3\}$.
Zhu and Bu~\cite{ZB} proved that if $G$ is a triangle-free planar graph, then $G$ is equitably $k$-choosable for $k\ge\max\{\Delta(G),8\}$ and if $G$ is a planar graph without $4$- and $5$-cycles, then $G$ is equitably $k$-choosable for $k\ge\max\{\Delta(G),7\}$.

It is worth mentioning that for every $k$, there are graphs that are equitably $k$-colorable but not equitably $k$-choosable, and vice versa. %$\textcolor{red}{(Examples?)}
One of the examples is the star $K_{1,2k}$.

\begin{claim}
    For each integer $k\ge 3$, the star $K_{1,2k}$ is equitably $k$-choosable but not equitably $k$-colorable.
\end{claim}

\begin{proof} Let $G=K_{1,2k}$ and let $v$ be the unique vertex of degree $2k$
in $G$.

        Assume that every vertex is assigned  a list of size $k$.
    Note that regardless of what color is used for $v$, there are at least $k-1$ colors to choose from for each neighbor of $v$.
    By the definition of equitable list coloring, each color can be used up to $\lceil(2k+1)/k\rceil=3$ times. 
    Since $3(k-1)\ge 2k$, we may color the neighbors of $v$ greedily so that no color is used more than $3$ times.

 If $K_{1,2k}$ had an
    equitable $k$-coloring, then  the color class of $v$ cannot be used 
 for any other vertex. Also in this case,  each of the remaining $k-1$ color classes would have at most two vertices. Thus $|G|\leq 1+2(k-1),$ a contradiction.
 %, we must use exactly one color $3$ times.
 %   Then we can color at most $3+2(k-2)=2k-1$ neighbors of $v$, so some vertex may not be colored.
\end{proof}

The next example of an equitably $k$-colorable but not equitably $k$-choosable graph is a bit more sophisticated.

\begin{ex}
    %For each integer $k\ge 3$, we construct a graph $G_k$ as follows. 
    { Consider any integer  $k\ge3$ and disjoint sets $V_1,U_1,\dots,V_k,U_k$ with $|V_i|=k-1$ and $|U_i|=k^3-2k^2+1$. Define the graph $G=G_k=(V,E)$  on $V:=\bigcup_{i=0}^k(V_i\cup U_i)$ by 
    \[E:=\bigcup_{i=0}^k\{vw:(v,w)\in ((V_0\times V_i)\cup (V_i\times U_i))\}.\]} 
   % $V=V(G_k)=\bigcup_{i=0}^k(V_i\cup U_i)$ 
   % where $|V_i|=k-1,|U_i|=k^3-2k^2+1$ and all $V_i,U_i$ are disjoint for $0\le i\le k$., so $|V|=k(k-1)^2(k+1)$.
    %For distinct $u,v\in V$, $uv\in E(G)$ if and only if $u\in V_0$, $v\in V_i$ or $u\in U_i,v\in V_i$ for some $0\le i\le k$.
\end{ex}
    \begin{claim}
        For each integer $k\ge 3$, $G_k$ is equitably $k$-colorable but not equitably $k$-choosable.
    \end{claim}
    \begin{proof}
   {
    Notice that
    \begin{enumerate}[label=(\alph*),nosep]
    \item $|G|=(k+1)|V_0\cup U_0|=(k+1)k(k-1)^2$, 
    \item $|N(V_0)|=k|V_1|+|U_0|=k^3-2k^2+1+k(k-1)=k^3-k^2-k+1=(k+1)(k-1)^2$, 
    \item $|V\smallsetminus N[V_0]|=(k-1)(k+1)(k-1)^2$,
    \item $V_0$ is a $(k-1)$-clique, and
    \item $N(V_0)$ and $V\smallsetminus N[V_0]$ are independent.
    \end{enumerate}
    Now $G$ is equitably $k$-colorable, since we can color   $V_0$ with $k-1$ distinct colors, use the same $k-1$ colors to equitably color $V\smallsetminus N(V_0)$, and color $N(V_0)$ with the remaining color.
    
    %After coloring the vertices in $V_0$, there is only one color $\alpha$ that can be used for any $w\in N(V_0)$. 
    %In total, 
    %Color the vertices in $U_1\cup\dots\cup U_k$ with colors other than $\alpha$ so that each color is used exactly $(k-1)^2(k+1)-1$ times.
    %There is no restriction on coloring those vertices, so we may use each color, including $\alpha$, exactly $(k-1)^2(k+1)$ times overall.}

    To see that $G$ is not equitably $k$-choosable, we assign a particular $k$-list $L$ to $V$ and show that $G$ does not have an equitable $L$-coloring.
    Set %$[k]=\{1,2,\dots,k\}$,
    $A:=[2k-1]\smallsetminus [k]$ and $V_i:=\{v_{i,j}:j\in [k-1]\}$ for $i\in [k]$.
    For $v\in V_0\cup U_0$, let $L(v)=[k]$.
    For each $v_{i,j}$, let $L(v_{i,j})=([k]-i)\cup \{k+j\}$.
    For $v\in U_i$ with $i\in [k]$, let $L(v)=A+i$.

    Since $V_0$ is a clique, we must use $k-1$ distinct colors from $[k]$ on $V_0$.
    Without loss of generality, we do not use the color $1$; then all vertices in $U_0$ must be colored $1$.
    By the definition of $L$, $v_{1,j}$ must get color $k+j$ for $j\in [k-1]$.
    Now every color in $A$ is used on $V_1$, so vertices in $U_1$ must get color $1$.
    Now there are $|U_0|+|U_1|=2k^3-4k^2+2$ vertices that must be colored with $1$, which already exceeds $|V|/k=(k-1)^2(k+1)$ when $k\ge 3$.}
    \end{proof}
%\end{ex}

%To unify these two notions, t
The authors~\cite{KKX} defined the \emph{strongly equitable} (SE) $L$-coloring which has the property that every SE $k$-choosable graph is also equitably $k$-choosable and equitably $k$-colorable. An $L$-coloring
 is SE if it has at most $|G|\mm k$ full classes, where $n\mm k$ is the unique
integer $m\in[k]$ with $n-m$ divisible by $k$. So, if $n$ is divisible by $k$
then $n\mm k=k$. A multigraph is SE $k$-choosable if it is SE $L$-colorable for
every $k$-list assignment $L$.
If a multigraph is SE $k$-choosable, then it is equitably $k$-colorable and equitably $k$-choosable.

Let $\mathcal{B}$ denote the class of graphs such that
\begin{equation}
\|X,Y\|\leq2|X\cup Y|-4~~\text{for all disjoint}~~X,Y\subseteq V(G)~~\text{with}~~|X\cup Y|\ge3.\label{eq:bibound}
\end{equation}
Observe that class $\mathcal{B}$ is topologically  much broader than the
class of planar graphs. For example, for any graph $H$, the graph $B_H$ obtained from $H$ by
subdividing each edge once is in $\mathcal{B}$. Therefore, every graph is a minor of a graph in $\mathcal{B}$. In
particular, the acyclic chromatic number of graphs in $\mathcal{B}$ can be arbitrarily large, while it is
at most 5 for any planar graph. Note that the class $\mathcal{Z}$
of graphs $G$ such that $||G[A]||\leq 4(|A|-2)$
for every $A\subseteq V(G)$ with $|A|\geq 3$ considered by
Zhang~\cite{1planar} also has these properties.

{ Adapting the method described in Section~\ref{maxd} to SE list coloring and analysing  sparsity of graphs in $\mathcal{B}$, we recently proved the following bound.}

%{\HK I THINK WE SHOULD  MENTION THAT THIS COMBINES OUR BASIC TECHNIQUE WITH SPARSE METHODS}
\begin{thm}[\cite{KKX}, 2024]\label{thm:eqlistb}
    Every graph $G\in \BB$ 
  is SE $k$-choosable,  if $k\geq\max\{9,\Delta(G)\}$. 
\end{thm}
Theorem~\ref{thm:eqlistb} combined with Theorem~\ref{thm:KK13}(a) implies the following fact.

\begin{cor}\label{co1}
  Conjecture~\ref{conj:list-HS} holds for every graph in $\mathcal{B}$, in particular, for every planar graph.  
\end{cor}

We think that it is natural to extend Conjectures~\ref{conj:list-HS},~\ref{conj:Orelist} and~\ref{conj:list-CLW}  to SE list coloring as follows.

\begin{conjecture}\label{conj:list-HS-SE}
    Every graph $G$ is equitably SE $k$-choosable for each $k> \Delta(G)$.
\end{conjecture}

\begin{conjecture}\label{conj:Orelist-SE}
  {If  $2k>\Theta(G)$, then $G$ is SE $k$-choosable.}
    %Every graph $G$ is SE $(1+0.5\Theta(G))$-choosable. 
\end{conjecture}

\begin{conjecture}\label{conj:list-CLW-SE}
    If $G$ is an $k$-colorable graph with $3\le\Delta(G)\le k$, then either $G$ is SE $k$-choosable or $k$ is odd and $G$ contains a $K_{k,k}$.
\end{conjecture}

\section{Coloring sparse graphs with fewer colors }\label{fewer}

One of the directions of study of equitable coloring and list equitable coloring is trying to color "sparse" graphs $G$ with fewer than $\Delta(G)$ colors.
By "sparse" one can mean "small average degree" or "small maximum average degree" or a topologically simple class of graphs, like outerplanar or planar graphs. 

Berge and Sterboul~\cite{BeS} 
studied equitable coloring of sparse graphs and more generally uniform hypergraphs. They considered
$m_0(n,r,k)$---the minimum number of edges in an $n$-vertex $r$-uniform hypergraph that has no equitable $k$-coloring. They have found the exact values of $m_0(n,r,k)$ for many triples $(n,r,k)$. In particular,
for graphs ($r=2$), they have found the exact values for all pairs $(n,k)$.
They proved that
$$m_0(n,2,k)=\left\{\begin{array}{ll}
\min\{{k+1\choose 2}, {n\choose 2}-{2(n-k)-1\choose 2}\}     &\mbox{when } k\leq n\leq 2k; \\
 \min\{{k+1\choose 2}, n-\lfloor\frac{n}{k}\rfloor
 +1  \}   &\mbox{when }  n\leq 2k.
\end{array}\right.
$$
The extremal example for small $k$  with respect to $n$ is simply $K_{k+1}$ that has no proper $k$-coloring at all together with $n-k-1$ isolated vertices. For "larger" $k$,  if $n\geq 2k$ then the extremal example is the star $K_{1,n-\lfloor\frac{n}{k}\rfloor
 +1}$ together with some isolated vertices, and if $k\leq n\leq 2k$, then it is the graph obtained from
 $K_n$ by deleting the edges of a $K_{2(n-k)-1}$.

Meyer~\cite{Mey} proved that every tree with maximum degree $\Delta$ has an equitable $k$-coloring
for $k = 1 + \lfloor \frac{\Delta}{2}\rfloor$ and it was observed later that the result holds for every $k \geq 1 + \lfloor \frac{\Delta}{2}\rfloor$. The extremal example is the star $K_{1,\Delta}$.

%The first result on equitable  $k$-coloring of sparse graphs
%whose maximum degree is not bounded in terms of $k$ is due to 
Recall that Theorem~\ref{BG} by Bollob\' as and Guy~\cite{BG83} gives for trees a bound that is not a function of maximum degree. For convenience, we restate it here.

\begin{thm}[Bollob\' as and Guy~\cite{BG83}, 1983]\label{BoG}
    A tree $T$ is equitably $3$-colorable if $|T|\ge 3\Delta(T)-8$ or $|T|=3\Delta(T)-10.$
\end{thm}
They also provided an algorithm for producing the coloring. 
An extremal example is any $n$-vertex tree containing a vertex $v$ of degree exactly $\frac{n+9}{3}$ and a perfect matching on the set of the nonneighbors of $v$.
This result was extended to all $k\geq 2$ and to all forests by Chen and Lih~\cite{CL94} and Miyata,  Tokunaga and Kaneko~\cite{MTK}.
Given a graph $G$ and a vertex $v\in V(G)$, let $\alpha_v=\alpha_v(G)$ denote the size of the maximum independent set in $G$ containing $v$. If a graph $G$ has an equitable $k$-coloring, then by definition, $\alpha_v(G)\geq \lfloor n/k\rfloor$ for every $v\in V(G)$.
Chen and Lih~\cite{CL94} and  independently Miyata,  Tokunaga and Kaneko~\cite{MTK} proved (and Chang~\cite{Ch08} gave a shorter and nicer proof) that this necessary condition is sufficient for forests.

\begin{thm}[Chen and Lih~\cite{CL94}, 1994, Miyata,  Tokunaga and Kaneko~\cite{MTK}, 1994]\label{CLCh}
    For a forest $T$ of order $n$ and integer $k\ge 3$, $T$ is equitably $k$-colorable if and only if $\alpha_v\ge\lfloor n/k\rfloor$ for every vertex $v\in T$.
\end{thm}

The "only if" part is obvious.
We present the proof of the "if" part due to  Chang from~\cite{Ch08}.

{\bf Proof.} Let $n=|V(G)|$.
We choose a bipartition $(A,B)$ of $V(G)$ so that $|A|=a\geq |B|=b$, and modulo this, $a$ is minimum. By this choice,
\begin{equation}\label{AB}
 \mbox{either $|A|-|B|\leq 1$ or $A$ does not contain isolated vertices.}   
\end{equation}
For $1\leq i\leq k,$ let $s_i=\lfloor\frac{n+i-1}{k}\rfloor.$
Choose the minimum $j$ such that $b\leq \sum_{i=1}^js_i$. If
$b=\sum_{i=1}^js_i$, then we partition $B$ into sets of sizes 
$s_1,\ldots,s_j$ and $A$ into sets of sizes 
$s_{j+1},\ldots,s_k$. So we may assume
$ \sum_{i=1}^{j-1}s_i<b< \sum_{i=1}^js_i$.

{\bf Case 1:} $j>1$. Take $S$ be the set of $s:=b-\sum_{i=1}^{j-1}s_i$ vertices of lowest degrees in $B$. By this choice,
$|N(S)|\leq \frac{s(n-1)}{b}$. Thus using   $a\geq b$ and $b-s\geq s_1$, we have
$$|S\cup (A-N(S))|>s+a-\frac{sn}{b}=\frac{(b-s)a}{b}\geq b-s\geq s_1.
$$
So, $|S\cup (A-N(S))|\geq 1+s_1\geq s_j$, and $S\cup (A-N(S))$ contains a subset $S'$ with $|S'|=s_j$ containing $S$. Then we can partition $A-S'$ into (independent) sets of sizes $s_{j+1},\ldots,s_k$.

\medskip
{\bf Case 2:} $j=1$. In this case, it is enough to prove that
\begin{equation}\label{Case2}
 \mbox{there are disjoint independent sets $I_1,I_2$ with $|I_1|\geq s_1$, $|I_2|\geq s_k$ and $B\subset I_1\cup I_2$,}   
\end{equation}
since in this case $V(G)-I_1-I_2$ is independent. Thus we assume that~\eqref{Case2} does not hold.

Since $k\geq 3$,~\eqref{AB} implies that $A$ does not contain isolated vertices. Let $L$ be the set of all leaves in $A$. Then
$|L|+2|A-L|\leq |E(G)|\leq n-1$ and so
$|L|\geq 2a-(n-1)=a-b+1
$. In particular,
\begin{equation}\label{LB}
|L|+b\geq a+1\geq n-b+1.
\end{equation}

Call a set $Q\subseteq B$ {\em good} if $|Q\cup (L-N(Q))|\geq s_k$.
Since~\eqref{Case2} does not hold, $B$ is not good. By~\eqref{LB}, $\emptyset$ is good. Let $Q$ be a largest good subset $Q$ of $B$. Since $B$ is not good, there is $v\in B-Q$. By the maximality of $Q$,
\begin{equation}\label{newv}
|Q+v|+|L-N(Q)-N(v)|\leq s_k-1.
\end{equation}
If $|(N(Q)\cap L)\cup (B-Q)|\geq s_1$, then we are done by~\eqref{Case2}. So, suppose $|N(Q)\cap L|+ |B-Q|\leq s_1-1$. Summing this with~\eqref{newv}, we get $|L|-|N(v)\cap L|+b+1\leq s_1+s_k-2$. Thus, using~\eqref{LB},
\begin{equation}\label{n(v)}
|N(v)\cap L|+(b-1)\geq |L|+b+1-s_1-s_k+2+(b-1)\geq n+2-s_1-s_k>s_k.
\end{equation}
Since $\alpha_v\geq s_1,$ there is an independent set $R$ containing $v$ with $|R|\geq s_1$. Choose such $R$ with the minimum $|R\cap B|$.
If $R\cap B=\{v\}$, then by~\eqref{n(v)}, there is an independent set $I_2\subseteq (B-v)\cup (N(v)\cap L)$ with $|I_2|=s_k$ disjoint from $R$, and we are done by~\eqref{Case2}.

Otherwise, $|R\cap B|\geq 2$, say $v'\in (R-v)\cap B$. If there is $x\in L-R$ with no neighbors in $R$, then the set $R-v'+x$ has fewer vertices in $B$ than $R$, contradicting the choice of $R$. So each $x\in L-R$ has a neighbor in $R$. Then the set $I_2:=(B\cup L)-R$ is independent, and by~\eqref{LB}, $|I_2|\geq n-b+1-s_1\geq s_k.$
Thus~\eqref{Case2} holds with $I_1=R$, a contradiction.
\qed

Chang~\cite{Ch08} asked about broader classes $\FF$ of graphs for which the conclusion of Theorem~\ref{CLCh} holds.

\begin{question}\label{qu1} Does the claim of Theorem~\ref{CLCh} hold for SE list coloring? In other words, is it true that for each $k\geq 3$,  every $n$-vertex forest $T$ is equitably SE $k$-choosable if and only if $\alpha_v\ge\lfloor n/k\rfloor$ for every vertex $v\in T$?  
\end{question}

\bigskip

When Yap and Zhang~\cite{YZ1} proved the CLW-Conjecture for outerplanar graphs, they also conjectured
 that every
outerplanar graph with maximum degree $\Delta\geq 3$ is equitably $k$-colorable for every $k\geq 1+\Delta/2$. The conjecture was proved in~\cite{Kout}.

\begin{question}\label{qu2} Does the claim of Theorem~\ref{CLCh} hold for all outerplanar graphs?
\end{question}

When Zhang and Wu~\cite{ZhW} proved the CLW-Conjecture and its list version for the broader class of
series-parallel graphs, they also conjectured
 that every series-parallel
 graph with maximum degree $\Delta\geq 3$ is equitably $k$-colorable for every integer $k\geq (\Delta+3)/2$.
 Kostochka and Nakprasit~\cite{KN03} proved a more general claim for $d$-degenerate graphs $G$ with $\Delta(G)\geq 27d$:

\begin{thm}[Kostochka and Nakprasit~\cite{KN03}, 2003]\label{KN03}
    Let $d\geq 2$, $\Delta\geq 27d$, $k \geq (d + \Delta + 1)/2$. Then every $d$-degenerate graph with
maximum degree at most $\Delta$ is equitably $k$-colorable.
\end{thm}
The restriction on $k$ is tight. Consider the graph $G(d,\Delta) = K_d +\overline{K}_{\Delta-d+1}$, obtained from the complete
graph $K_d$ by adding $\Delta-d+1$ vertices each of which is adjacent to vertices
of our $K_d$ and only to them. This graph is $d$-degenerate, has maximum degree $\Delta$ and needs $\left\lceil\frac{\Delta+d+1}{2}\right\rceil$ colors for its equitable coloring because some $d$   color classes must be of size $1$.

In particular, Theorem~\ref{KN03} confirms the conjecture by Zhang and Wu for
$\Delta\geq 54$. The restriction $\Delta\geq 27d$ in Theorem~\ref{KN03} maybe is not needed, and Kostochka and Nakprasit~\cite{KN03} conjectured that the claim of Theorem~\ref{KN03} holds for all $2\leq d\leq \Delta$. Partially supporting this conjecture is the result of
 Verclos and Sereni~\cite{VS17} who proved the conjecture by Zhang and Wu in full.

Pemmaraju~\cite{Pe} proved an analogue  of Theorem~\ref{BoG} for outerplanar graphs.

\begin{thm}[Pemmaraju~\cite{Pe}, 2001]\label{Pe1}
   A connected outerplanar $n$-vertex graph with  maximum 
degree at most $n/6$ has an equitable $6$-coloring.
\end{thm}

 Kostochka, Nakprasit and Pemmaraju~\cite{KNP} generalized this as follows.

\begin{thm}[Kostochka, Nakprasit and Pemmaraju~\cite{KNP}, 2005]\label{knp}
For $d, n \geq  1$, every $d$-degenerate $n$-vertex graph $G$ with 
$\Delta(G)\leq n/15$ 
is equitably $k$-colorable for each $k \geq 16d$.
\end{thm}

Again, maybe the restrictions $\Delta(G)\leq n/15$ and $k \geq 16d$ can be weakened (but cannot be dropped completely).

In a similar spirit, Wu and
Wang \cite{WW} proved in 2008 that every planar graph $G$ with minimum degree $\delta(G)\geq2$
and girth $g(G)\geq14$ is equitably $k$-colorable for each $k\geq4$. Moreover,
if the girth of $G$ is at least $26$ then $G$ is equitably $3$-colorable. Luo,
Sereni, Stephens and Yu \cite{LSSY} strengthened this result with the following
theorem.
\begin{thm}[Luo, Sereni, Stephens and Yu \cite{LSSY}, 2010]
\label{thm:LuoEtAl}Every planar graph $G$ with $\delta(G)\geq2$ and $g(G)\geq10$
is equitably $k$-colorable for each $k\geq4$, and if the $g(G)\geq14$ then $G$
is equitably $3$-colorable. 
\end{thm}

Luo \emph{et al.} \cite{LSSY} 
%%%AK added citation
mention that their proofs yield the same bounds when they replace
planarity with appropriate bounds on maximal average degree: for graphs $G$ with
$\mad(G)\leq7/3$ when $k=3$ and for graphs with $\mad(G)\leq5/2$ when $k\geq4$,
but they still need the restrictions on the girth and minimum degree.

%The authors~\cite{KKX2} sharpened, generalized and extended to SE list colorings Theorem~\ref{thm:LuoEtAl} for $k\leq 4$ as follows. 
{ Theorem~\ref{thm:LuoEtAl} was sharpened, generalized and extended to SE list colorings as follows.}

\begin{thm}[\cite{KKX2}, 2024]
\label{thm:main0}Let $G$ be a graph with $\delta(G)\geq2$. If 
\begin{equation}\label{3col}
6\|G[A]\|\leq 7|A|+2\quad\mbox{for each nonempty $A\subseteq V(G)$} %HK \| not ||
\end{equation}
%$G$ is $(\frac{7}{6},\frac{1}{3})$-sparse
then $G$ SE $3$-choosable, and if 
\begin{equation}\label{4col}
4\|G[A]\|\leq 5|A|+2\quad\mbox{for each nonempty $A\subseteq V(G)$}
\end{equation}
%$G$ is $(\frac{5}{4},\frac{1}{2})$-sparse
  then
$G$ is SE $4$-choosable. Moreover, both results are sharp.
\end{thm}

%Recall that a girth of a graph $g(G)$ is the length of its shortest %cycle.

\begin{figure}
\begin{center}
\begin{tikzpicture}[scale =1] 
\def \vt{circle (1.5pt) [fill]} \def \mg{cyan!50!red!20!white} \def \mc{red!20!yellow!20!white} \def \mpg{cyan!20!green!20!white} \foreach \j in {-10,35,80,125,170} { \draw (\j:2cm)--(\j+20:2cm)--(0:0cm)--cycle; 
\draw (\j:2cm) \vt; \draw (\j:1cm) \vt; \draw (\j+20:2cm) \vt; \draw (\j+20:1cm) \vt;} 
\draw (0,0) \vt; 
\draw (0,0)--(.5,-.5)--(.5,-1)--(.25,-1.5)--(-.25,-1.5)--(-.5,-1)--(-.5,-.5)--cycle; 
\draw (.5,-.5) \vt; 
\draw (.5,-1) \vt;
\draw (-.5,-1) \vt; 
\draw (-.5,-.5) \vt; 
\draw (.25,-1.5) \vt;
\draw (-.25,-1.5) \vt;
\draw (0,0) \vt; 
\end{tikzpicture}
\qquad
\begin{tikzpicture}[scale =1]
\def \vt{circle (1.5pt) [fill]} \def \mg{cyan!50!red!20!white} \def \mc{red!20!yellow!20!white}  \def \mpg{cyan!20!green!20!white} \foreach \j in {20,60,100,140} { \draw (\j:2cm)--(\j+20:2cm)--(0:0cm)--cycle; 
\draw (\j:2cm) \vt; \draw (\j+20:2cm) \vt; } 
\draw (0,0) \vt; 

\draw (0,0)--(.5,-.5)--(.5,-1)--(-.5,-1)--(-.5,-.5)--cycle; 
\draw (.5,-.5) \vt; 
\draw (.5,-1) \vt;
\draw (-.5,-1) \vt; 
\draw (-.5,-.5) \vt; 
\draw (0,0) \vt;
\end{tikzpicture}
\end{center}

\caption{\protect\label{fig:k=00003D3} Graph  on the left of the picture has $27$ vertices, $32$ edges and cannot be equitably $3$-colored.
Graph  on the right of the picture has $13$ vertices, $17$ edges and cannot be equitably $4$-colored.
}\label{fig3}
\end{figure}
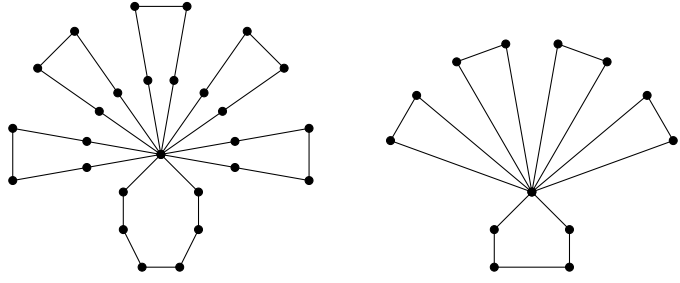

Sharpness examples for Theorem~\ref{thm:main0} are { shown} in Fig.~\ref{fig3}. It is interesting that the bounds for ordinary equitable coloring, list equitable coloring and SE list coloring are the same.
We think that we can prove an exact bound for SE list coloring with $5$ colors, but are not sure what would be the exact bound for more colors.

Krivelevich and Patkos~\cite{KrP} studied equitable colorings of random graphs in the model $G(n,p)$. They showed that if
$n^{-1/5+\epsilon} < p < 0.99$ for some $\epsilon>0$, then almost surely the minimum $k$ such that $G\in G(n,p)$ has an equitable $k$-coloring is asymptotically the chromatic number of $G$. They also proved that the maximum $k$ such that if
$C/n < p < 0.99$ for a large $C$, then almost surely the maximum $k$ such that $G\in G(n,p)$ has no equitable $k$-coloring is asymptotically at most twice the chromatic number of $G$.

\section{Variations of equitable coloring}

In this section we discuss some variants of equitable coloring.
Recall that in equitable list coloring, sizes of color classes could differ as some colors may appear only  in few lists.
Given a graph $G$ and a $k$-list assignment $L$ for $G$, for each color $c$, denote  by $\eta(c)$  the number of vertices $v\in V(G)$ such that $c\in L(v)$.
A \textit{proportional $L$-coloring} of $G$ is a proper $L$-coloring that uses each color $c$ exactly $\lfloor \eta(c)/k\rfloor$ or $\lceil \eta(c)/k\rceil$ times.
A graph $G$ is \textit{proportionally $k$-choosable} if for every $k$-list assignment $L$, $G$ has a proportional $L$-coloring.
The \textit{proportional choice number} $\chi_{pc}(G)$ is the smallest $k$ such that $G$ is proportionally $k$-choosable.

This notion of proportional coloring was proposed by Kaul, Mudrock, Pelsmajer and Reiniger in~\cite{KMPR}, and they proved the following upper bound on $\chi_{pc}(G)$ for all graphs $G$.

\begin{thm}[Kaul, Mudrock, Pelsmajer and Reiniger~\cite{KMPR}, 2018]
    For every graph $G$, $$\chi_{pc}(G)\le \Delta(G)+\left\lceil\frac{|V(G)|}{2}\right\rceil.$$
\end{thm}

They also asked if the analog of Theorem~\ref{thm:H-Sz} holds for proportional coloring. 

\begin{question}[Kaul, Mudrock, Pelsmajer and Reiniger~\cite{KMPR}, 2018]
    For any graph $G$, is $G$ proportionally $k$-choosable whenever $k\ge \Delta(G)+1$?
\end{question}

Kaul \emph{et al.} \cite{KMPR} mention that the analog of Conjecture~\ref{conj:CLW} does not hold for proportional coloring. 
{A} %The 
counterexample is a union of $k$ disjoint $K_{1,k}$ for each integer $k$.

Coloring of a graph could be viewed as a partition of the vertices into disjoint independent sets.
A \textit{$k$-partition} of a graph is a partition of the vertex set into $k$ disjoint sets.
There have been many studies on partitioning the vertex set into sets other than independent sets, for example, sets that induce graphs with
 {bounded} %certain
 degeneracy or maximum degree.
Similarly to equitable coloring, an\textit{ equitable $k$-partition} is a partition into $k$ disjoint sets such that the sizes of these sets differ by at most $1$.

Kostochka, Nakprasit and Pemmaraju  proved the following result on partitioning a $d$-degenerate graph into $(d-1)$-degenerate graphs equitably.

\begin{thm}[Kostochka, Nakprasit and Pemmaraju~\cite{KNP}, 2005]
    For $k\ge 3,d\ge 2$, every $d$-degenerate graph has an equitable $k$-partition into $(d-1)$-degenerate graphs.
\end{thm}

By this result, each planar graph can be partitioned into $81$ forests since every planar graph is $5$-degenerate.
Esperet, Lemoine, and Maffray found a better upper bound on the number  $k$ of forests such that each planar graph could be partitioned into $k$ forests, and asked whether an even better upper bound could be achieved.

\begin{thm}[Esperet, Lemoine, and Maffray~\cite{ELM}, 2015]
    For each $k\ge 4$, every planar graph has an equitable $k$-partition into forests.
\end{thm}

\begin{question}[Esperet, Lemoine, and Maffray~\cite{ELM}, 2015]
    Does every planar graph have an equitable $3$-partition into forests?
\end{question}

Relaxing the degeneracy of the  sets in a partition, Kim, Oum and Zhang derived the following results.

\begin{thm}[Kim, Oum and Zhang~\cite{KOZ}, 2022]
Every planar graph has 
\begin{enumerate}[label=(\alph*),nosep]
    \item an equitable $2$-partition into $3$-degenerate graphs,
    \item  an equitable $3$-partition into $2$-degenerate graphs,
    \item an equitable $3$-partition into two forests and {a} %one
     graph.
\end{enumerate}
\end{thm}

If a vertex partition $(V_1,\ldots,V_k)$ of a graph $G$ is such that the maximum degree of
$G[V_i]$
is at most  $d$ for each $1\leq i\leq k$, then such partition is a \textit{$d$-defective coloring} of $G$.
%If a graph $G$ has an equitable $k$-partition that is $1$-defective, then $G$ is \textit{equitable defective $k$-colorable}.
 Williams, Vandenbussche and Yu  proved the following result on equitable defective coloring that is an extension of Theorem~\ref{thm:LuoEtAl}.

\begin{thm}[Williams, Vandenbussche and Yu~\cite{WVY}, 2012]
    For every $k\ge 3$ each planar graph $G$ with $\delta(G)\geq2$ and $g(G)\geq10$
has an equitable $k$-partition that is a $1$-defective coloring.
\end{thm}

\bibliographystyle{plain}
\bibliography{Ref}

\end{document}